\documentclass[11pt]{article}

\usepackage {amsmath}
\usepackage {amsfonts}
\usepackage {amsthm}
\usepackage[parfill]{parskip}
\usepackage {amssymb}
\usepackage {mathrsfs}

\usepackage {fullpage}

\usepackage {amscd}
\usepackage[all,cmtip]{xy}

\title{Modular Categories associated to Unipotent Groups}
\author{Tanmay Deshpande}
\date{}

\newtheorem {thm} {Theorem} [section]
\newtheorem {prop} [thm] {Proposition}
\newtheorem {conj} [thm] {Conjecture}
\newtheorem {lem} [thm] {Lemma}
\newtheorem {cor} [thm] {Corollary}
\theoremstyle{definition}
\newtheorem {defn} [thm] {Definition}
\newtheorem {example} [thm]  {Example}
\theoremstyle{remark}
\newtheorem {rk} [thm]  {Remark}

\renewcommand {\bar} {\overline}
\newcommand{\bpf}{\begin{proof}}
\newcommand{\epf}{\end{proof}}

\newcommand{\f}{\mathbb}

\newcommand{\h}{\hbox}
\newcommand{\ind}{\h{ind}}

\renewcommand{\L}{\mathcal{L}}

\newcommand{\Fp} {\mathbb{F}_p}

\newcommand{\tK} {\widetilde{K}}

\newcommand{\Q} {\mathbb{Q}}

\newcommand{\W} {\mathbb{W}}

\newcommand{\F} {\mathbb{F}}
\newcommand{\normal} {\triangleleft}
\newcommand{\tN}{\widetilde{N}}

\newcommand{\FPdim} {\hbox{FPdim}}
\renewcommand{\k} {\mathtt{k}}
\newcommand{\K} {\mathcal{K}}
\newcommand{\M} {\mathcal{M}}
\newcommand{\N} {\mathscr{N}}

\newcommand{\D} {\mathcal{D}}
\newcommand{\C} {\mathcal{C}}

\newcommand{\Hom} {\hbox{Hom}}
\renewcommand{\Vec}{\mbox{Vec}}
\newcommand{\Aut}{\mbox{Aut}}
\newcommand{\Out}{\mbox{Out}}
\newcommand{\Inn}{\mbox{Inn}}

\newcommand{\Qlcl} {\overline{\mathbb{Q}}_l}
\newcommand{\beq}{\begin{equation}}
\newcommand{\eeq}{\end{equation}}
\newcommand{\bthm}{\begin {thm}}
\newcommand{\ethm}{\end {thm}}
\newcommand{\bprop}{\begin {prop}}
\newcommand{\eprop}{\end {prop}}
\newcommand{\bcor}{\begin {cor}}
\newcommand{\ecor}{\end {cor}}
\newcommand{\blem}{\begin{lem}}
\newcommand{\elem}{\end{lem}}
\newcommand{\bdefn}{\begin{defn}}
\newcommand{\edefn}{\end{defn}}
\newcommand{\brk}{\begin{rk}}
\newcommand{\erk}{\end{rk}}

\newcommand{\G} {\mathbb{G}}
\newcommand{\noin}{\noindent}

\newcommand{\tg} {\tilde{\gamma}}
\newcommand{\ta} {\tilde{a}}
\newcommand{\tb} {\tilde{b}}
\newcommand{\tab} {\widetilde{ab}}
\renewcommand{\t} {\widetilde}
\newcommand{\tc} {\tilde{c}}
\newcommand{\td} {\tilde{d}}
\newcommand{\tga} {\tilde{\gamma_1}}
\newcommand{\tgab} {\widetilde{\gamma_1\gb}}
\newcommand{\tgb} {\tilde{\gamma_2}}
\newcommand{\g} {{\gamma}}
\renewcommand{\l} {{\lambda}}

\newcommand{\ga} {{\gamma_1}}
\newcommand{\gb} {{\gamma_2}}

\newcommand{\cpu}{\mathfrak{cpu}}

\newcommand{\cpuc}{\mathfrak{cpu}^\circ}

\newcommand{\Meg}{\M_{G,e}}
\newcommand{\tM}{\widetilde{\M}}
\newcommand{\tH}{\widetilde{H}}
\newcommand{\bit}{\begin{itemize}}
\newcommand{\eit}{\end{itemize}}
\mathchardef\hy = "2D

\newcommand{\onto}{\twoheadrightarrow}

\newcommand{\E}{\mathcal{E}}
\newcommand{\Z}{\mathcal{Z}}
\renewcommand{\Vec}{\mbox{Vec}}
\newcommand{\bconj}{\begin{conj}}
\newcommand{\econj}{\end{conj}}
\newcommand{\Gcpuc}{\Gamma\hy\cpuc}
\newcommand{\Gcpu}{\Gamma\hy\cpu}
\newcommand{\Gap}{\G_{a,perf}}
\newcommand{\Id}{\mbox{Id}}

\newcommand{\uuPic}{\underline{\underline{\hbox{Pic}}}}
\newcommand{\uPic}{\underline{\hbox{Pic}}}
\newcommand{\Pic}{\hbox{Pic}}

\newcommand{\uEqBr}{{\underline{\hbox{EqBr}}}}
\newcommand{\EqBr}{{\hbox{EqBr}}}

\begin{document}
\maketitle

\begin{abstract} 
\noin Let $G$ be a unipotent algebraic group over an algebraically closed field $\k$ of characteristic $p>0$ and let $l\neq p$ be another prime. Let $e$ be a minimal idempotent in $\D_G(G)$, the $\Qlcl$-linear triangulated braided monoidal category of $G$-equivariant (for the conjugation action) $\Qlcl$-complexes on $G$ under convolution (with compact support) of complexes. Then, by a construction due to Boyarchenko and Drinfeld,  we can associate to $G$ and $e$ a modular category $\M_{G,e}$. In this paper, we prove that the modular categories that arise in this way from unipotent groups are precisely those in the class $\mathfrak{C}_p^{\pm}$.
\end{abstract}

\section{Introduction} \label{introduction}
Let $G$ be a unipotent group over an algebraically closed field $\k$ of characteristic $p>0$. We fix the field $\k$ and all algebraic groups and schemes we consider will be assumed to be over $\k$. By an algebraic group, we mean a smooth group scheme of finite type over the field $\k$. A unipotent group is an algebraic group (over $\k$) that is isomorphic to a closed subgroup of the algebraic group $UL_n(\k)$ of unipotent upper triangular matrices of size $n$ for some positive integer $n$. 

Let us also fix a prime number $l\neq p$. For a $\k$-scheme $X$ we have the $\Qlcl$-linear triangulated symmetric monoidal category $\D(X):=D^b_c(X,\Qlcl)$  of constructible $\Qlcl$-complexes on $X$ under tensor product of complexes. For any integer $n$, we have the degree shift by $n$ functor $(\cdot)[n]:\D(X)\to\D(X)$ as well as the $n$-th Tate twist functor $(\cdot)(n):\D(X)\to\D(X)$. If we have an action $\alpha:G\times X\to X$ of the unipotent group $G$ on $X$, we have the $G$-equivariant $\Qlcl$-linear triangulated cateogry $\D_G(X)$. We recall (from \cite[\S 1.3]{BD08}) that an object of $\D_G(X)$ is a pair $(M,\phi)$ where $M\in\D(X)$ and $\phi:\alpha^*M\to p_2^*M$ is an isomorphism (known as the $G$-equivariant structure) satisfying certain compatibility conditions, where $p_2:G\times X\to X$ is the second projection. We have the degree shift and Tate twist functors on $\D_G(X)$ as well.

For the unipotent group $G$, the $\Qlcl$-linear triangulated category $\D(G)$ is also monoidal under convolution with compact support. We recall (from \cite[\S 1.4]{BD08}) that for $M,N\in\D(G)$, we define the convolution with compact support as $M\ast N:=\mu_!(M\boxtimes N)$ where $\mu:G\times G\to G$ is the group multiplication and $\mu_!$ denotes (derived) pushforward with compact support. Here the external tensor product $M\boxtimes N:=p_1^*M\otimes p_2^*N$ where $p_i:G\times G\to G$ for $i=1,2$ are the two projections. We also have the $\Qlcl$-linear triangulated {\it braided monoidal category} $\D_G(G)$ of $G$-equivariant complexes for the conjugation action of $G$ on itself.  For $M\in\D(G)$ and $(N,\phi)\in \D_G(G)$, we have braiding isomorphisms $\beta_{M,N}:M\ast N\to N\ast M$ defined in \cite[Defn. A.43]{BD08} which provide the braided structure on $\D_G(G)$. The category $\D_G(G)$ is equipped with a twist $\theta$, an automorphism of the identity functor on $\D_G(G)$ obtained using the $G$-equivariant structure of the objects of $\D_G(G)$ (see \cite[\S4.9]{B},\cite[\S5.5]{BD06} or \cite[\S1.5.4, Appendix A.5]{BD08}). This gives $\D_G(G)$ the structure of a ribbon $\mathfrak{r}$-category. We refer to \cite[\S5.5]{BD06}, \cite[\S1 and Appendix A.5]{BD08} for more on the structure of the categories $\D(G)$ and $\D_G(G)$.  

Let $e\in\D_G(G)$ be a minimal idempotent. This means that $0\neq e\cong e\ast e$ and for any idempotent $e'\in\D_G(G)$, $e\ast e'$ is either isomorphic to $e$ or 0. We refer to \cite[Defn. 1.11, also Thm. 1.49]{BD08} for the precise definition and some subtler points. Associated to the minimal idempotent $e$, we have the $\Qlcl$-linear triangulated Hecke subcategory $e\D_G(G)$, which is a braided monoidal category with unit object $e$. We recall from {\it loc. cit.} that $e\D_G(G)$ is the full subcategory of $\D_G(G)$ consisting of objects $M\in \D_G(G)$ such that $M\cong e\ast N$ for some $N\in \D_G(G)$ (or equivalently, such that $M\cong e\ast M$).  In \cite{BD08}, Boyarchenko and Drinfeld define an $\f{L}$-packet of character sheaves on $G$ and a modular category $\M_{G,e}\subset e\D_G(G)$ corresponding to each such minimal idempotent $e$. We recall this in Theorem \ref{BDmain0} below.

\brk
All monoidal categories that we consider in this paper will be $\Qlcl$-linear. In particular all the fusion categories that we talk about are assumed to be over the field $\Qlcl$. We recall that a fusion category over $\Qlcl$ is a $\Qlcl$-linear semisimple abelian rigid tensor category with finitely many simple objects (up to isomorphism), finite dimensional $\Hom$ spaces and with unit object being simple. A spherical structure on a fusion category is a certain additional structure (see \cite[\S2.4.3]{DGNO}) which, in particular, allows us to define categorical dimensions $d(X)\in\Qlcl$ of objects $X$ in the fusion category. A modular category is a non-degenerate braided spherical category. We refer to \cite{ENO02} and \cite{DGNO} for more on the notions of fusion categories, spherical structures, braided fusion categories, ribbon structures and modular categories.
\erk

\bthm\label{BDmain0}{(\cite{BD08})}
Let $G$ be a unipotent group over the field $\k$. Let $e\in \D_G(G)$ be a minimal idempotent. Then the following hold:
\bit
\item [(i)] Let $\M_{G,e}^{perv}\subset e\D_G(G)$ denote the full subcategory consisting of objects of $e\D_G(G)$ whose underlying $\Qlcl$-complex is a perverse sheaf on $G$. Then $\M_{G,e}^{perv}$ is a semisimple abelian category with finitely many simple objects. The (isomorphism classes of) simple objects of $\M_{G,e}^{perv}$ are said to form the $\f{L}$-packet of character sheaves associated to the minimal idempotent $e$.

\item [(ii)] There exists an integer $n_e\in \{0,1,\cdots,\dim(G)\}$ such that $e\in \M_{G,e}^{perv}[n_e]$. Moreover, the shifted subcategory $\M_{G,e}:=\M_{G,e}^{perv}[n_e]\subset e\D_G(G)$ is closed under convolution and the twist $\theta$ in $\D_G(G)$ induces the structure of a modular category on $\M_{G,e}$. Moreover, we have $D^b(\M_{G,e})\cong e\D_G(G).$

\eit
\ethm

\bdefn
The number $d_e:=\frac{\dim(G)-n_e}{2}\in \frac{1}{2}\f{Z}$ is said to be the functional dimension of the minimal idempotent $e$.
\edefn

\noin Thus if we have a unipotent group $G$, and a minimal idempotent $e\in \D_G(G)$, we have the associated modular category $\Meg$. In this paper, we characterize all modular categories that arise in this way from unipotent groups, namely these modular categories are precisely those in the class $\mathfrak{C}_p^{\pm}$. We state the technical definition and some properties of this class of modular categories in Appendix \ref{B}. (See also \cite{DGNO2}.) We briefly recall the definition of the class $\mathfrak{C}_p^+$ below.

A pointed fusion category is one in which all simple objects are invertible. Consequently, the set of isomorphism classes of simple objects is a finite group $\Gamma$. All pointed fusion categories (over $\Qlcl$) are of the form $\Vec_\Gamma^{\omega}$, the category of finite dimensional $\Gamma$-graded $\Qlcl$-vector spaces with usual tensor product, but with the associativity constraint defined by a 3-cocycle
$\omega\in H^3(\Gamma, \Qlcl^*)$. A pointed fusion category has a unique positive spherical structure (see \cite[Ex. 2.26]{DGNO}) and with this structure all simple objects are of categorical dimension 1. The Drinfeld center $\Z(\D)$ of such a pointed (positive spherical) fusion category $\D$ is a modular category, which can also be realized as the category of representations of the so-called twisted Drinfeld double of a finite group. A modular category $\C$ is said to be in the class $\mathfrak{C}_p^+$ if we have an equivalence of modular categories $\C\cong \Z(\Vec_{\Gamma}^\omega)$, where $\Gamma$ is a {\it finite $p$-group} and $\Vec_{\Gamma}^\omega$ is equipped with the positive spherical structure. (See also Proposition \ref{dgno2}.) The definition of the class $\mathfrak{C}_p^-$ is more technical and we provide it in Appendix \ref{B}.

\noin The goal of this paper is to prove the following result conjectured by Drinfeld in \cite{BD06}:

\bthm\label{conj1}
(i) Let $G$ be a unipotent group. Let $e\in \D_G(G)$ be a minimal idempotent with functional dimension $d_e$ ($\in \frac{1}{2}\f{Z}$). Then if $d_e\in \mathbb{Z}$, $\M_{G,e}\in \mathfrak{C}_p^{+}$ and if $d_e\in \frac{1}{2}+\mathbb{Z}$,  $\M_{G,e}\in \mathfrak{C}_p^{-}$.\\
(ii) If $\C\in \mathfrak{C}_p^{\pm}$, then there exists a \emph{connected} unipotent group $G$ with a minimal idempotent $e\in \D_G(G)$ such that $\M_{G,e}\cong \C$ as modular categories.
\ethm

\noin We prove Theorem \ref{conj1}(i) in \S\ref{pfi}. In \cite{BD08}, Boyarchenko and Drinfeld prove that every minimal idempotent $e\in \D_G(G)$ can be induced from a {\it Heisenberg idempotent} (Definition \ref{defhi}) $e'\in \D_{G'}(G')$ satisfying a certain {\it geometric Mackey condition}, where $G'$ is some subgroup of $G$. They also prove that $d_e=d_{e'}+\dim(G/G')$ and that there is an equivalence $\Meg\cong \M_{G',e'}$ of modular categories. In fact, we prove in \S\ref{redhi} that we may choose $e'$ to be a {\it special} (Definition \ref{defhi}) Heisenberg idempotent. Hence we reduce Theorem \ref{conj1}(i) to the case where $e$ is a  special  Heisenberg idempotent. In \S\ref{pfhi} we complete the proof of Theorem \ref{conj1}(i) by proving it in the case of special Heisenberg idempotents.

In \S\ref{pfii}, we prove Theorem \ref{conj1}(ii). Let $\C\in \mathfrak{C}_p^{\pm}$. In \S\ref{dtoc}, we show that it is enough to find a possibly disconnected unipotent group $G$ with a Heisenberg idempotent $e\in \D_G(G)$ such that $\C\cong \Meg$. We show that if $G$ is a (possibly disconnected) unipotent group with a Heisenberg idempotent $e$, then there exists a {\it connected} unipotent group $U'$ and a minimal idempotent $f'\in \D_{U'}(U')$ such that $\Meg\cong \M_{U',f'}$. Finally, we complete the proof of Theorem \ref{conj1}(ii) for the case $\C\in \mathfrak{C}_p^{+}$ in \S\ref{cp+} and the case $\C\in \mathfrak{C}_p^{-}$ in \S\ref{cp-}

\subsection*{Acknowledgments}
I would like to thank V. Drinfeld for introducing me to this subject, for the many useful discussions and for suggesting corrections and improvements.

\section{Proof of Theorem \ref{conj1}(i)} \label{pfi}
In this section, we prove the first part of the main theorem. We begin by reducing to the case of Heisenberg idempotents and further to the case of  special Heisenberg idempotents.

\subsection{Reduction to special Heisenberg idempotents} \label{redhi}
\noin Let us recall the notions of admissible pairs and Heisenberg idempotents defined in \cite{BD08}. A multiplicative local system on a {\it connected} unipotent group $G$ is a nonzero $\Qlcl$-local system $\L$ on $G$ equipped with an isomorphism $\mu^*\L\cong \L\boxtimes \L$. It follows that $\L$ must be of rank 1, and that it has a natural $G$-equivariant structure and we can consider it as an object of $\D_G(G)$. According to \cite[Lem. 7.3]{B}, we can equivalently think of a multiplicative local system on $G$ as a central extension of $G$ by the abstract group $\Q_p/\f{Z}_p$ if we choose an identification of $\Q_p/\f{Z}_p$ with $\mu_{p^\infty}$, the group of $p$-power-th roots of unity in $\Qlcl^*$.

\brk\label{multloc}
In particular, any central extension $0\to C\to \tilde{G}\to G\to 0$ of $G$ by an abelian group $C$ along with a character $C\to \mu_{p^\infty}$ determines a multiplicative local system on $G$.  
\erk

It will often be necessary to work with perfect schemes. For a scheme $X$ over $\k$, we denote by $X_{perf}$ its perfectization. We have $X(\k)=X_{perf}(\k)$. For a perfect connected commutative unipotent group $U$ over $\k$, we let $U^*$ denote its Serre dual which parametrizes the multiplicative $\Qlcl$-local systems on $U$. $U^*$ is also a perfect connected commutative unipotent group. (See \cite[Appendix A.5]{B}, \cite[\S1]{Be} for details.)

\bdefn\label{defap}(\cite[Defn. 1.30]{BD08})
Let $G$ be a unipotent group. Let $(N,\L)$ be a pair consisting of a connected subgroup $N\subset G$ and a multiplicative local system $\L$ on $N$. We say that the pair $(N,\L)$ is admissible if the following conditions hold:
\bit
\item[(i)] Let $G'=\N_G(N,\L)$, the normalizer of the pair $(N,\L)$ in $G$ (in particular $ N \normal G'$ and $\L$ is $G'$-equivariant, see also \cite[Defn. 1.29]{BD08}) and let $G'^\circ$ denote its neutral connected component. Then $G'^{\circ}/N$ is commutative.
\item[(ii)] The group morphism $\phi_\L:(G'^{\circ}/N)_{perf} \to (G'^{\circ}/N)_{perf}^*$ that is defined in this situation (see Remark \ref{ssi} below) is an isogeny.
\item[(iii)] {\it(Geometric Mackey condition)} For every $g\in G(\k)-G'(\k)$, we have $$\L|_{(N\cap N^g)^\circ}\ncong\L^g|_{(N\cap N^g)^\circ},$$ where $N^g=g^{-1}Ng$ and $\L^g$ is the multiplicative local system on $N^g$ obtained from $\L$ by transport of structure. 
\eit
\edefn

\bdefn\label{defsap}
Let $(N,\L)$ be an admissible pair as above. If we further have that $\dim(G'/N)\leq 1$, we say that $(N,\L)$ is a {\it special} admissible pair.

\edefn

\brk
For an admissible pair $(N,\L)$, let $e':=\L[2\dim N](\dim N)\in \D_{G'}(G')$. Then $e'\in \D_{G'}(G')$ is in fact a minimal idempotent on the subgroup $G'\subset G$ (see \cite[Prop. 8.2]{B}, \cite{De}). 
\erk

\brk\label{ssi}
In view of statement (i) of the definition, we have the commutator map $G'^{\circ}\times G'^\circ\to N$. The pullback of $\L$ via the commutator induces a skew-symmetric biextension of $(G'^\circ/N)_{perf}\times (G'^\circ/N)_{perf}$. (See \cite[Appendix A.6]{B} for more on biextensions.) This defines the map $\phi_\L$ of (ii). More explicitly, for $g\in G'^\circ$, let $c_g:G'^\circ\to N$ be the {\it commutator with $g$} map: $G'^\circ\ni h \mapsto hgh^{-1}g^{-1}$. Then we can check that $c_g^*\L$ is a multiplicative local system on $G'^\circ$ whose restriction to the normal subgroup $N$ is trivial, whence we can think of it as a multiplicative local system on $G'^\circ/N$. Explicitly, the map $\phi_\L: (G'^{\circ}/N)_{perf} \to (G'^{\circ}/N)_{perf}^*$ is defined by $g\mapsto c_g^*\L$. Often, we say that the map $\phi_\L$ itself is skew-symmetric. Also the conjugation action of $G'$ on $G'^\circ$ induces an action of $\pi_0(G')$ on $G'^{\circ}/N$. Since $\L$ is $G'$-equivariant, it follows that $\phi_\L$ is $\pi_0(G')$-equivariant. See \cite[\S 3.3]{BD08} for more.
\erk

\bdefn\label{defhi}
If $(N,\L)$ is an admissible pair for $G$ such that $G'=G$, we say that $(N,\L)$ is a Heisenberg admissible pair\footnote{Note that in this situation, condition (iii) of the definition is vacuous.} for $G$. In this case the minimal idempotent $e:=\L[2\dim N](\dim N)\in \D_G(G)$ is said to be a Heisenberg idempotent. Further, if we have $\dim(G/N)\leq 1$ we say that the idempotent is a special Heisenberg idempotent.
\edefn

\brk
It is easy to see that if $e\in\D_G(G)$ is a Heisenberg idempotent as above, then $n_e=\dim(N)$. This is because $\L[\dim(N)]$ is a perverse sheaf. Hence the functional dimension $d_e=\frac{1}{2}\dim(G^\circ/N)$.
\erk

\begin{example}\label{ex}
Let us consider an instructive example of the {\it fake Heisenberg group} $H$. As a variety, $H=\G_a\times \G_a$ and we define multiplication by $(x,a)\cdot(y,b)=(x+y,a+b+xy^p)$. Note that $(x,a)^{-1}=(-x,x^{p+1}-a).$ The center of $H$ is $N:=0\times \G_a$ and $H/N\cong \G_a$ is commutative. Let $\L$ be the multiplicative local system on $N$ corresponding to the Artin-Schreier central extension $0\to \F_p\to \G_a\xrightarrow{x\mapsto x^p-x} \G_a\to 0$ (and the choice of a primitive $p$-th root of unity $\zeta\in\Qlcl^*$). Since $N$ is central, it is clear that $\N_H(N,\L)=H$. The commutator with $(y,b)$ map $c_{(y,b)}:H\to N\cong \G_a$ defined in Remark \ref{ssi} is given by $c_{(y,b)}(x,a)=xy^p-x^py=:c_y(x)$ (since in this example the commutator map is independent of $a,b$). Let us pass to the perfectization, $\Gap$ of $\G_a$. We can identify $\Gap$ with its Serre dual $\Gap^*$. This identification $\psi$ can be described as follows (see \cite[\S1.1]{Be} for details):\\
For $t\in \Gap$, we have the homothety $t:\Gap\to\Gap$. Then we define $\psi(t)=t^*\L$, the pullback (considered as an element of $\Gap^*$) of the Artin-Schreier multiplicative local system by the homothety $t$.\\
We see that in this example, the map $\phi_\L:\Gap\to \Gap^*$ is defined by $y\mapsto c_y^*\L$. If we identify $\Gap^*$ with $\Gap$ as above, we can check that this map is given by $y\mapsto y^p-y^{{1}/{p}}$. This map is an isogeny with kernel $\F_{p^2}$. Hence $(N,\L)$ is a Heisenberg admissible pair on $H$. The object $e:=\L[2](1)$ is the corresponding Heisenberg idempotent and $d_e=\frac{1}{2}$.
\end{example}

\noin In \cite[\S 1.12]{BD08}, the induction (with compact support) functor $\hbox{ind}_{G'}^G:\D_{G'}(G')\to \D_G(G)$ is defined for closed subgroups $G'\subset G$. In the following theorem, Boyarchenko and Drinfeld prove that every minimal idempotent in $\D_G(G)$ comes from an admissible pair and in particular from a Heisenberg idempotent on a subgroup by induction.

\bthm\label{BDmain}
\cite{BD08}
(i) Let $(N,\L)$ be an admissible pair for a unipotent group $G$ and let $e'\in \D_{G'}(G')$ be the corresponding Heisenberg idempotent on $G'$ as defined above. Then $e:=\ind_{G'}^Ge'\in \D_G(G)$ is a minimal idempotent. The functional dimensions of $e$ and $e'$ are related as follows:$$d_e=d_{e'}+\dim(G/G').$$\\
(ii) In the situation of (i), the induction functor induces an equivalence of modular categories $\M_{G',e'}\cong \M_{G,e}$.\\
(iii) Every minimal idempotent $e\in \D_G(G)$ comes from an admissible pair by the procedure described in (i). Hence every minimal idempotent $e\in \D_G(G)$ comes from induction from a Heisenberg idempotent $e'$ on some subgroup $G'$.
\ethm

Let us now strengthen this last statement and prove that every minimal idempotent can in fact be induced from a special Heisenberg idempotent. Since we know that every minimal idempotent can be induced from a Heisenberg idempotent, it suffices to prove that every Heisenberg idempotent can be induced from a special Heisenberg idempotent. For ease of notation, henceforth in this section let us work in the setting of {\it perfect} unipotent groups, namely all groups considered are assumed to be perfect unipotent. We will use results proved in Appendix \ref{A}. 

Let us suppose we have a Heisenberg admissible pair $(N,\L)$ on a perfect unipotent group $G$. Let $e=\L[2\dim N](\dim N)\in\D_G(G)$ be the corresponding Heisenberg idempotent. Let $H=G^\circ$ and let $\Gamma=G/H=\pi_0(G)$. Since $G$ is unipotent, $\Gamma$ is a $p$-group. Note that we have the $\Gamma$-equivariant skew-symmetric isogeny $\phi_\L:H/N\to (H/N)^*$. 

\bprop
Let $W$ be a connected subgroup of $H$ containing $N$ such that $W/N$ is a $\Gamma$-invariant isotropic subgroup of $H/N$. Then we can extend $\L$ to a multiplicative local system $\L'$ on $W$. The pair $(W,\L')$ is admissible for $G$.
\eprop
\bpf
Since $W/N$ is isotropic, the pullback of $\L$ by the commutator map $c_W:W\times W\to N$ is trivial. Hence by \cite[Prop. 7.7]{B} we can extend $\L$ to a multiplicative local system $\L'$ on $W$. 

The $\Gamma$-invariance of $W/N$ implies that $W\normal G$.  Let $W'$ be the connected subgroup of $H$ such that $W'/N = ((W/N)^\perp)^\circ$. Then we see that $\N_G(W,\L')^\circ = W'$ and $W'/W$ is commutative. Note that the commutator map $c_{W'}:W'\times W'\to W$ factors as $c_{W'}:W'\times W'\to N\subset W$. It follows that the skew-symmetric biextension $\phi_{\L'}:W'/W\to (W'/W)^*$ defined using the pair $(W,\L')$ is the same as the skew-symmetric biextension of $W'/W = ((W/N)^{\perp})^\circ/(W/N)$ induced from the isotropic subspace $W/N$ for skew-symmetric biextension $\phi_\L:H/N \to (H/N)_{perf}^*.$ It follows that $\phi_{\L'}$ is an isogeny. Since $W\normal G$, the geometric Mackey condition is immediate. Hence $(W,\L')$ is an admissible pair for $G$.
\epf

\bcor
Let $G$ be a unipotent group. Every minimal idempotent $e\in\D_G(G)$ comes from a special admissible pair.
\ecor
\bpf
Let $(N,\L)$ be an admissible pair for a minimal idempotent $e$. We assume without loss of generality that the normalizer of the pair $(N,\L)$ equals $G$. (Else we can replace $G$ by the normalizer of the pair.) We will now continue to use previous notations. Let us now choose a connected subgroup $W\subset H$ containing $N$ such that $W/N$ is a $\Gamma$-invariant maximal isotropic subgroup of $H/N$. The existence of $W$ is guaranteed by Corollary \ref{mgiis}. Then we can extend $\L$ to a multiplicative local system $\L'$ on $W$. In this case $\N_G(W,\L')^{\circ}=W'$ as defined above. Since $W/N$ is maximal isotropic, $\dim(W'/W)\leq 1$ by Corollary \ref{mgiis} and hence $(W,\L')$ is a special admissible pair for the minimal idempotent $e$.
\epf

\subsection{Analysis of Heisenberg idempotents}\label{anhi}
Let $G$ be a unipotent group. Let $(N,\L)$ be a Heisenberg admissible pair for $G$ and let $e=\L[2\dim N](\dim N)\in\D_G(G)$ be the corresponding Heisenberg idempotent. Let $H=G^\circ$ and let $\Gamma=G/H$. We see that in this case $n_e=\dim (N)$ and hence $d_e=\frac{\dim(H/N)}{2}$. We have the $\Gamma$-equivariant skew-symmetric isogeny $\phi_\L:(H/N)_{perf}\to (H/N)_{perf}^*$. Let $K_\L=\h{ker}\phi_\L$. As defined in \cite[Appendix A.10]{B} or \cite[\S 2.3]{Da} we have a quadratic form (see Definition \ref{defmg}) $\theta:K_\L\to \Qlcl^*$ associated to the skew-symmetric isogeny $\phi_\L$ and we have a metric group $(K_\L,\theta)$. We also have an action of $\Gamma$ on this metric group.  According to Appendix \ref{mcassmg}, we have the pointed positive spherical modular category $\M(K_\L,\theta)$ associated to the metric group $\M(K_\L,\theta)$. We see in Theorem \ref{Demain} below that $\M_{H,e}\cong \M(K_\L,\theta)$. 

Consider the full subcategory $e\D_H(G)\subset \D_H(G)$, the $\Qlcl$-linear triangulated monoidal category of $H$-equivariant complexes. According to \cite[\S 2.4]{De}, we have a grading 
$$e\D_H(G)=\bigoplus\limits_{H\tg\in G/H}e\D_H(H\tg),$$ 
a monoidal action of $\Gamma$ on $e\D_H(G)$ and a {\it crossed braiding} which give $e\D_H(G)$ the structure of a braided $\Gamma$-crossed category (see \cite[\S 4.4.3]{DGNO} for definition, see also Appendix \ref{bcce}) with identity component $e\D_H(H)$ and unit object $e$. The equivariantization (see \cite[\S 4.4.4]{DGNO}) $e\D_H(G)^\Gamma$ is a braided monoidal category equivalent to the Hecke subcategory $e\D_G(G)\subset \D_G(G)$.

Let $\tM_{G,e}$ denote the full subcategory of $e\D_H(G)$ consisting of perverse sheaves shifted by $\dim N$. Note that by definition $\M_{H,e}\subset \tM_{G,e}$ is the full subcategory consisting of shifted (by $\dim(N)$) perverse sheaves supported on $H\subset G$. The following theorem is proved in \cite{De}.

\bthm\label{Demain}
(\cite{De})
(i) $\tM_{G,e}$ is closed under convolution and we have a faithfully graded braided $\Gamma$-crossed structure  on $\tM_{G,e}$ with trivial component $\M_{H,e}$.\\ 
(ii) $\M_{H,e}$ is a pointed modular category equivalent to the modular category $\M(K_\L,\theta)$ associated to the metric group $(K_\L,\theta)$. \\
(iii) The equivariantization ${\tM_{G,e}}^{\Gamma}$ is a modular category and $\M_{G,e}\cong {\tM_{G,e}}^{\Gamma}$.
\ethm

\brk\label{action}
By statements (i) and (ii) above, we have a braided action of $\Gamma$ on $\M(K_\L,\theta)$ coming from the braided $\Gamma$-crossed structure of $\tM_{G,e}$. This induces an action of $\Gamma$ on the metric group $(K_\L,\theta)$ which is the same as the action which we had defined before.
\erk

\brk\label{heisid}
We note that whenever we have a minimal idempotent $e$ on a unipotent group $G$, we have the modular category $\M_{G,e}$. Moreover, if $e$ is a Heisenberg idempotent as described above, we also have a braided $\Gamma$-crossed category $\tM_{G,e}$ whose identity component is the pointed modular category $\M_{H,e}\cong \M(K_\L,\theta)$ and whose equivariantization $(\tM_{G,e})^\Gamma$ is equivalent to $\M_{G,e}$.
\erk

\begin{example}\label{ex1}
Let $H$ be the fake Heisenberg group with the Heisenberg admissible pair $(N,\L)$ and corresponding Heisenberg idempotent $e$ defined in Example \ref{ex}. We saw that in this example, the skew-symmetric isogeny 
$$\phi_\L:\Gap\to\Gap^*\cong\Gap$$
is given by $y\mapsto y^p-y^{1/p}$ and that the kernel $K_\L=\F_{p^2}$. In this case, (see \cite[Appendix A]{Da}) the metric group associated to the skew-symmetric isogeny is the anisotropic metric group $(\F_{p^2},\zeta^{Nm(\cdot)})$ where $Nm:\F_{p^2}\to \F_p$ is the norm map. Hence in this example, the modular category $\M_{H,e}$ is equivalent to $\M_p^{anis}:=\M(\F_{p^2},\zeta^{Nm(\cdot)})$.
\end{example}

\brk\label{spherical}
The category $\tM_{G,e}$ is a rigid monoidal category. According to \cite[\S 2.3]{De} the duality functor is given by $\tM_{G,e}\ni X\mapsto X^\vee=\f{D}\iota^*X[2\dim N](\dim N)$ where  $\f{D}$ denotes Verdier duality and $\iota:G\to G$ is the inversion map.  We have natural isomorphisms $X\cong (X^\vee)^\vee$ which give us a canonical pivotal (in fact spherical) structure on $\tM_{G,e}$. Hence ${\tM_{G,e}}^{\Gamma}\subset e\D_G(G)$ also has a spherical structure which is the same as the one obtained using the twist $\theta$ on $\D_G(G)$.
\erk

In the remainder of this subsection, we study certain invariants of the modular category $\M_{G,e}$ in case of Heisenberg and special Heisenberg idempotents. We refer to \cite[\S 2.4.1]{DGNO} for the notion of Frobenius-Perron dimension and \cite[\S 6]{DGNO} for the notions of Gauss sums and multiplicative central charge.

We prove the following proposition for all Heisenberg idempotents, and later restrict to the class of special Heisenberg idempotents.  Recall that for a Heisenberg idempotent $e\in\D_G(G)$, $d_e=\frac{1}{2}\dim(H/N)$. We continue to use the assumptions and notations introduced at the beginning of this subsection.

\bprop\label{prop1}
Let $e\in\D_G(G)$ be a Heisenberg idempotent.\\
(i) The Gauss sums $\tau^{\pm}(\Meg)$ are equal to $(-1)^{2d_e}p^{k}\cdot |\Gamma|$ for some integer $k\geq 0$. Hence the multiplicative central charge of $\Meg$ equals $(-1)^{2d_e}$.\\
(ii) The Frobenius-Perron dimension of $\Meg$ equals $|\Gamma|^2\cdot |K_\L|$.
\eprop
\bpf
Since $\M_{H,e}$ is equivalent to the modular category corresponding to the metric group $(K_\L,\theta)$, by \cite[Prop. 5]{Da} we have $\tau^{\pm}(\M_{H, e,})=\tau^\pm(K_\L,
\theta)=(-1)^{2d_e}p^k$ for some integer $k\geq 0$. 
Suppose $(M,\phi)\in \M_{H,e}^{\Gamma}\subset \Meg$, where $M\in \M_{H,e}$ and $\phi$ is the equivariant structure. The twist $\theta_{(M,\phi)}$ in $\Meg$ is the same as the twist $\theta_M$ in $\M_{H,e}$ as an automorphism of $M$. Let $\mathcal{E}$ be the full subcategory of $\M_{H,e}^{\Gamma}$ consisting of objects $(M,\phi)$ such that $M$ is isomorphic to a direct sum of copies of the unit object $e$. Hence it follows that the twist when restricted to $\mathcal{E}$ is trivial, or in other words, $\mathcal{E}$ is an isotropic subcategory of $\Meg$. Moreover, it is clear that $\mathcal{E}\cong \mbox{Rep}(\Gamma)$. The braided $\Gamma$-crossed category $\tM_{G,e}$ is precisely the corresponding de-equivariantization and the identity component $\M_{H,e}$ is the corresponding fiber category. (See \cite[\S 4.4.7-8]{DGNO}.) Hence by \cite[Thm. 6.16]{DGNO}, $\tau^\pm(\Meg)=\tau^\pm(\M_{H,e})\cdot|\Gamma|$. Statement (i) now follows, since  $\tau^\pm(\M_{H,e})=(-1)^{2d_e}p^k$. To prove (ii), note that $\mbox{FPdim}(\tM_{G,e})=|\Gamma|\cdot \mbox{FPdim}(\M_{H,e})=|\Gamma|\cdot|K_\L|$ by \cite[Prop. 2.21]{DGNO} and 
that $\mbox{FPdim}(\tM_{G,e}^\Gamma)=|\Gamma|\cdot\mbox{FPdim}(\tM_{G,e})$ by \cite[4.26]{DGNO}. 
\epf

\brk
Since $(K_\L,\theta)$ is a metric group coming from a skew-symmetric biextension, it follows from $\cite{Da}$ that $|K_\L|$ is an even power of $p$. Moreover, $\Gamma$ is a $p$-group. Hence we see that $\h{FPdim}(\M_{G,e})=p^{2k}$ for some $k\in \f{Z}^+$. We also see that the multiplicative central charge of $\Meg$ is 1 if $\dim(H/N)$ is even and $-1$ if $\dim(H/N)$ is odd.
\erk

\brk 
Suppose the Heisenberg admissible pair $(N,\L)$ is special, i.e. $\dim(H/N)\leq 1$. Hence $H/N=0$ or $\f{G}_a$. Hence the induced action of the $p$-group $\Gamma$ on $H/N$ is trivial. 
\erk

Now suppose that we have any admissible pair such that the action of $\Gamma$ on $H/N$ is trivial. (By the above remark, every special admissible pair is of this type.) The action of $\Gamma$ on the metric group $(K_\L,\theta)$ is clearly trivial too. Since $\Gamma$ acts trivially on $H/N$, the commutator map from $H\times G$ in fact maps to $N$. We will now use some of the tools developed in \cite{De}. For  $g\in G$ we have the commutator map $c_g:H\to N$ defined $h\mapsto hgh^{-1}g^{-1}$. We recall from \cite[\S4.1]{De} that in this situation, the objects of $e\D_H(G)$ are supported on $\widetilde{K}:=\{g\in G|c_g^*\L\cong \Qlcl\}\supset N$. The objects of $e\D_H(H)$ are supported on  
$K:=\widetilde{K}\cap H$. By Remark \ref{ssi}, we see that $K/N=K_\L=\ker(\phi_\L)$. We can readily modify the proofs of Proposition 4.5, 4.6 from \cite{De}  and prove the following:

\bprop\label{cdo}
Let $e\in\D_G(G)$ be a Heisenberg idempotent as above such that the action of $\Gamma$ on $H/N$ is trivial. Then we have the following:\\
(i) $\tM_{G,e}$ is a \emph{central} (see Definition \ref{central}) braided $\Gamma$-crossed category with identity component $\M(K_\L,\theta)$.\\
(ii) Let $k\in \tK\subset G$.  Let $e^k$ denote the right translate of $e$ by $k$. Then $e^k\in \tM_{G,e}\subset e\D_H(G)$. The isomorphism class of $e^k$ only depends on the coset $Nk\in \tK/N$ and the $e^k, k\in \tK$ are all the simple objects of the braided $\Gamma$-crossed category $\tM_{G,e}$. \\
(iii) For $k_1,k_2\in \tK$ we have an isomorphism $e^{k_1}\ast e^{k_2}\cong e^{k_1k_2}$. Hence $\tM_{G,e}$ is pointed with the group of isomorphism classes of simple objects being $\tK/N$ which fits into a central extension $$0\to K_\L=K/N \subset\widetilde{K}/N\to \Gamma\to 0.$$ 
\eprop

\bpf
Since the action of $\Gamma$ on the metric group $(K_\L,\theta)$ obtained from the braided $\Gamma$-crossed category $\tM_{G,e}$ is trivial (see also Remark \ref{action}), $\tM_{G,e}$ is a central braided crossed category by Definition \ref{central}.

For $g\in G$, let $\delta_g$ denote the skyscraper sheaf supported at $g\in G$ whose stalk at $g$ is $\Qlcl$. By definition, the right translate $e^g=e\ast \delta_g$. Note that since $e\in\D_G(G)$, we have the braiding isomorphism $\delta_g\ast e\cong e\ast\delta_g$ between the left and right translates of $e$ by $g\in G$. For $k\in \tK$, $c_k^*\L\cong \Qlcl$ and hence we can define an $H$-equivariant structure on the translate $e^k=\L^k[2\dim N](\dim N)$. Hence $e^k\in \tM_{G,e}\subset e\D_H(G)$. For $n\in N$, a trivialization of the stalk $\L_n$ of $\L$ at $n$ defines an isomorphism $e\to e^n$ and also isomorphisms $e^k\to e^{nk}$ for each $k\in \tK$. Hence the isomorphism class of $e^k$ only depends on the coset $Nk$. It is clear that we have an isomorphism $e^{k_1}\ast e^{k_2}\cong e^{k_1k_2}$. In particular each $e^k\in\tM_{G,e}$ is an invertible simple object. From this, or using \cite[Prop. 3.10, 3.11]{De} we see that the $e^k$ for $k\in\tK$ are all the simple objects of $\tM_{G,e}$.
\epf

\bprop\label{cdpi}
The categorical dimensions of all the simple objects of $\tM_{G,e}$ are 1, or in other words, the spherical structure on the pointed category $\tM_{G,e}$ is positive (and of course integral). Hence the categorical dimensions of all simple objects of $\Meg\cong \tM_{G,e}^{\Gamma}$ are positive integers.
\eprop

\bpf
We must prove that the categorical dimension $d(e^k)=1$ for all $k\in \tK$. Note that $e^k$ is supported on $Nk$ and that its dual $(e^k)^\vee=\f{D}\iota^*(e^k)[2\dim N](\dim N)$ (see \cite[\S 2.3]{De} and Remark \ref{spherical} above) is supported on $k^{-1}N$. Here $\f{D}$ denotes Verdier duality and $\iota:G\to G$ is the inversion map. Let us identify $Nk$ and $k^{-1}N$ with $N$ by the evident maps. Under this identification $\iota:Nk\to k^{-1}N$ gets identified with $\iota:N\to N$, $\mu:Nk\times k^{-1}N\to N$ gets identified with the multiplication for $N$ (hence this identification is compatible with convolutions), $e^k$ gets identified with $e$ and Verdier duality is also compatible. By definition, the categorical dimension $d(e^k)$ equals the composition $e\to e^k\ast (e^k)^\vee\to (e^k)^{\vee\vee}\ast {(e^k)^\vee}\to e$. Once we make the identifications as above, we see that $d(e^k)=d(e)$ and $e$ is the unit object of $\tM_{G,e}$.  Hence categorical dimensions of all simple objects of $\tM_{G,e}$ are 1. Moreover, the categorical dimensions of all simple objects of $\Meg\cong\tM_{G,e}^\Gamma$ are positive integral multiples of those of $\tM_{G,e}$ (see \cite[Prop. 4.26]{DGNO}). Hence the categorical dimensions of all simple objects of $\Meg$ must also be positive integers.
\epf

\subsection{Completion of the proof of Theorem \ref{conj1}(i)}\label{pfhi}

Let $e\in\D_G(G)$ be a minimal idempotent. Then there exists a special Heisenberg idempotent $e'$ (coming from a special admissible pair $(N,\L)$ on $G$)  on a subgroup $G' \subset G$ satisfying the geometric Mackey condition such that $e\cong ind_{G'}^{G}e'$. By  Theorems \ref{BDmain}(ii) and \ref{Demain}(iii), we have $\M_{G,e}\cong\M_{G',e'}\cong \tM_{G',e'}^{\pi_0(G')}$. Since $e'$ is special Heisenberg idempotent on $G'$, we can apply Propositions \ref{cdo} and \ref{cdpi} and conclude that $\tM_{G',e'}$ is a central braided $\pi_0(G')$-crossed category, the spherical structure on $\tM_{G',e'}$ is positive and the identity component is equivalent to $\M(K_\L,\theta)$, where $(K_\L,\theta)$ is the metric group coming from the admissible pair $(N,\L)$. We recall that $K_\L=\ker(\phi_\L)\subset (G'^\circ/N)_{perf}$. Since $G'$ is unipotent, $\pi_0(G')$ is a $p$-group. By Theorem \ref{BDmain}(i), $d_e=\frac{\dim(G'^\circ/N)(=0 \hbox{ or }1)}{2}+\dim{G/G'}$. If $d_e\in \f{Z}$, $G'^\circ/N=0$ and then $\M(K_\L,\theta)\cong \Vec$.  Hence in this case $\M_{G,e}\in \mathfrak{C}_p^{+}$ by Definition \ref{defcppm}(i). If $d_e\in \frac{1}{2}+\f{Z}$, then $G'^\circ/N\cong \G_a$ whence the metric group $(K_\L,\theta)$ is Witt-equivalent to the anisotropic metric group (see Appendix \ref{mgassssbe}) $(\F_{p^2},\zeta^{Nm(\cdot)})$ by Theorem \ref{swarny}. Hence in this case $\M_{G,e}\in\mathfrak{C}_p^{-}$ by Definition \ref{defcppm}(ii). This completes the proof of Theorem \ref{conj1}(i).

\section{Proof of Theorem \ref{conj1}(ii)} \label{pfii}
Let $\C\in \mathfrak{C}_p^{\pm}$. In order to prove Theorem \ref{conj1}(ii), we first show that it is enough to find a possibly disconnected unipotent group $G$ with a Heisenberg idempotent $e$ such that $\C\cong\M_{G,e}$.

\subsection{Passing from a disconnected group to a connected group}\label{dtoc}
Let $G$ be a possibly disconnected unipotent group with a  Heisenberg admissible pair $(N,\L)$. As before, let $H=G^0$. Let $e\in\D_{G}(G)$ denote the corresponding Heisenberg idempotent. We will prove that there exists a connected unipotent group $U'$ and a minimal idempotent $f'\in \D_{U'}(U')$ such that the modular categories $\M_{G,e}$ and $\M_{U',f'}$ are equivalent.

\bprop\label{prep}
Let $V$ be a unipotent group with an action of $G$ by automorphisms. Let $\K$ be a $G$-equivariant multiplicative local system on $V$. Then we have the following:
\bit
\item[(a)] $(V\rtimes N, \K\boxtimes \L)$ is a Heisenberg admissible pair for $V\rtimes G$.
\item[(b)] We have an equivalence $e\D_{G}(G)\cong (f\boxtimes e)\D_{(V\rtimes G)}(V\rtimes G)$, where $f=\K[2\dim V](\dim V)$, given by $M\mapsto f\boxtimes M$. Here $f\boxtimes e$ is the Heisenberg idempotent corresponding to the admissible pair above.
\eit
\eprop
\bpf
Here $\K\boxtimes \L=p_1^*\K\otimes p_2^*\L$ is the external tensor product where $p_1, p_2$ are the two projections from $V\rtimes N$ to $V$ and $N$ respectively. Since $\K$ is $G$-equivariant multiplicative local system (and hence $N$-equivariant) on $V$, it is easy to check that $\K\boxtimes \L$ is a multiplicative local system on $V\rtimes N$ and that in fact it is also $V\rtimes G$-equivariant. We also check that the corresponding morphism $\left((V\rtimes H)/(V\rtimes N)\right)_{{perf}}\to\left((V\rtimes H)/(V\rtimes N)\right)_{{perf}}^*$ can be identified with the map $(H/N)_{{perf}}\to(H/N)_{{perf}}^*$ corresponding to the admissible pair $(N,\L)$. For part (b), we use the fact that $f$ is $G$-equivariant.
\epf


Let us now proceed to construct the connected unipotent group $U'$. First, we embed $G$ into a connected unipotent group $U$. Recall the well known theorem (due to Chevalley) that given an algebraic group along with a closed subgroup, we can find a representation of the group on a vector space and a line in the vector space such that the subgroup is precisely the stabilizer of the line. (See for example \cite[\S11.2]{H}.) In case the subgroup is unipotent, it can be characterized as the stabilizer of a point in the vector space (any nonzero point on the chosen line). Hence we can find a vector space $V$ with a $U$-action and a multiplicative local system $\K$ on $V$ (in other words $\K\in V^*$, which can be thought of as the vector space we start with) such that $G$ is precisely the stabilizer of $\K\in V^*$ under the $U$-action. Now let $U'=V\rtimes U$. Let $N'=V\rtimes N$ and $\L'=\K\boxtimes\L$. We will now prove that the pair $(N',\L')$ is admissible for $U'$, and that its normalizer is $G':=V\rtimes G$.

\bprop
\bit
\item[(a)] For $u\in U'-G'$ we have $$\L'|_{(N'\cap  { }^{u}N')^0}\ncong { }^u\L'|_{(N'\cap  { }^{u}N')^0}.$$
\item[(b)] $G'$ is the normalizer of the pair $(N',\L')$ and it is an admissible pair for $U'$. 
\item[(c)] Let $f'\in \D_{U'}(U')$ be the corresponding minimal idempotent. Then we have an equivalence $e\D_G(G)\cong f'\D_{U'}(U')$.
\eit
\eprop
\bpf
(a) Note that since $V$ is a connected normal subgroup of $U'$, we have $V\subset (N'\cap {}^uN')^0$. We have $\L'|_V\cong \K$ and that ${}^u\L'|_V\cong {}^u\K$. We have chosen $\K$ such that the stabilizer of $\K\in V^*$ under the action of $U$ is precisely $G$. Hence we have that ${}^u\K\ncong \K$ for each $u\in U'-G'$. Hence (a) follows.

\noin(b) It is clear that $G'$ is contained in the normalizer of $(N',\L')$. But from (a) we see that it is precisely the normalizer. This combined with (a) and Prop. \ref{prep} implies that $(N',\L')$ is an admissible pair for $U'$.

\noin(c) Note that we have $f'=\hbox{ind}_{G'}^{U'}(f\boxtimes e)$, where $f$ is an is Prop. \ref{prep}. Hence (c) follows from Prop. \ref{prep}(b) and Theorem \ref{BDmain}(ii).
\epf

\subsection{The case $\C\in \mathfrak{C}_p^{+}$}\label{cp+}
Let $\C\in \mathfrak{C}_p^{+}$. By Proposition \ref{dgno2}, we have a realization $\C\cong \Z(\Vec_\Gamma^{\omega})$ equipped with positive spherical structure for some finite $p$-group $\Gamma$ and an $\omega\in H^3(\Gamma,\Qlcl^*)$. Let us choose a representative $3$-cocycle $\omega:\Gamma^3\to \Qlcl^{\times}$.  Since $\Gamma$ is a finite $p$-group, we may assume that $\omega$ takes values in a finite subgroup $\mu_{p^n}\subset\Qlcl^{\times}$ of $p^n$-th roots of unity for some $n$. Let $\W$ denote the ring scheme of Witt vectors of length $n$. With its additive structure, $\W$ is a connected commutative unipotent group such that $p^n\cdot \W=0$. We have the Artin-Schreier sequence $0\to C:=\f{Z}/{p^n\f{Z}}\to \W\to \W\to 0$. We fix an identification $\mu_{p^n}\cong C$.

We will now construct a unipotent group $G$ along with a Heisenberg idempotent $e$ on it such that (see Remark \ref{heisid}) $\tM_{G,e}=\Vec_{\Gamma}^{\omega}$ as braided crossed categories (see proof of Proposition \ref{dgno2}). Hence as required we will have $\M_{G,e}\cong (\Vec_{\Gamma}^{\omega})^\Gamma\cong\Z(\Vec_{\Gamma}^{\omega})\cong \C$.  Consider the group ring $\tH:=\W(\Gamma)=Maps(\Gamma,\W)$. This is a connected commutative unipotent group (isomorphic to $\W^{|\Gamma|}$) with an action of $\Gamma$ by left translations. We embed $C\hookrightarrow \W\stackrel{\Delta}{\hookrightarrow}\tH$ via the diagonal. Note that under this embedding $\Gamma$ acts trivially on $C\subset \tH$. Let $H:=\tH/C$ as a $\Gamma$-module. Hence we have an exact sequence
\beq
0\to C \to \tH \to H\to 0
\eeq
of $\Gamma$-modules. It is a $\Gamma$-equivariant central extension of $H$ by $C\subset \Qlcl^*$. Let $\L$ denote the corresponding $\Gamma$-equivariant multiplicative local system on $H$. Note that since $\tH$ is the group algebra, $H^i(\Gamma, \tH)=0$ for $i>0$ by Shapiro's lemma. Hence from the long exact sequence of cohomology groups associated to the above short exact sequence, we get an isomorphism $\delta: H^2(\Gamma, H)\to H^3(\Gamma, C)$. Hence there exists an $f\in H^2(\Gamma, H)$ thats maps to $\omega$ (considered now as an element of $H^3(\Gamma,C)$) and we get a corresponding extension
\beq
0\to H\to G\to \Gamma\to 0.
\eeq

Note that $(H,\L)$ is a special Heisenberg admissible pair for $G$. Let $e$ be the corresponding idempotent.

For $\g\in\Gamma$ let us choose lifts $\tg\in G$ such that $\tga\tgb=f(\ga,\gb)\tgab$ where we view $f$ as a 2-cocycle.  We see that by Proposition \ref{cdo}, the translates $e^{\tg}$ are all the simple objects of $\tM_{G,e}$ (up to isomorphism). For $a,b\in \Gamma$, let us fix a trivialization of the stalk $\L_{f(a,b)}$. This determines an isomorphism $\eta_{a,b}:e\to e^{f(a,b)}$. For $a,b,c\in \Gamma$ we have ${ }^af(b,c)+f(a,bc)=f(a,b)+f(ab,c)=:h.$ From the $\Gamma$-equivariant multiplicative local system structure of $\L$ and the chosen trivializations $\eta_{\cdot,\cdot}$ we get two trivializations of $\L_h$. These two trivializations differ by $\delta f=\omega$ by definition of the boundary map. The simple objects of $\tM_{G,e}$ convolve as $$e^{\ta}\ast e^{\tb}\stackrel{\cong}{\longrightarrow} e^{\ta\tb}=e^{f(a,b)\tab}\cong e^{f(a,b)}\ast e^{\tab}\stackrel{\eta_{a,b}^{-1}}{\longrightarrow}e^{\tab}.$$
For $a,b,c\in \Gamma$, we have identifications $e^{\ta\tb\tc}\cong (e^{\ta}\ast e^{\tb})\ast e^{\tc}\cong e^{f(a,b)}\ast e^{f(ab,c)}\ast e^{\widetilde{abc}}$ and $e^{\ta\tb\tc}\cong e^{\ta}\ast (e^{\tb}\ast e^{\tc})\cong e^{{ }^a f(b,c)}\ast e^{f(a,bc)}\ast e^{\widetilde{abc}}$. Each of the previous two compositions are further identified with $e^{\widetilde{abc}}$ using the $\Gamma$-equivariant structure of $e$ and the chosen isomorphisms $\eta_{\cdot,\cdot}$. These two identifications differ by $\omega$. Hence the associativity constraint of $\tM_{G,e}$ corresponds to $\omega\in H^3(\Gamma,\Qlcl^*)$ i.e. $\tM_{G,e}\cong \Vec_{\Gamma}^\omega$ and by Proposition \ref{cdpi}, the spherical structure on $\tM_{G,e}$ is positive. Hence by Theorem \ref{Demain}(iii) after equivariantization, $\M_{G,e}\cong (\Vec_\Gamma^\omega)^\Gamma\cong$\footnote{See proof of Proposition \ref{dgno2}.} $\Z(\Vec_{\Gamma}^\omega)\cong \C$ as modular categories.

This (along with \S\ref{dtoc}) proves Theorem \ref{conj1}(ii) in the case when $\C\in \mathfrak{C}_p^{+}$.

\subsection{The case $\C\in \mathfrak{C}_p^{-}$}\label{cp-}
Let $\C\in \mathfrak{C}_p^{-}$. By Proposition \ref{dgno2}, there exists a faithfully graded\footnote{All braided crossed categories that we consider abstractly, as well as those that arise from other structures in this paper are faithfully graded even though sometimes we may not state this explicitly.} central (Definition \ref{central}) braided $\Gamma$-crossed category $\D$ (equipped with positive spherical structure) with identity component $\M_p^{anis}:=\M(\F_{p^2},\zeta^{Nm(\cdot)})$ such that $\C\cong \D^\Gamma$  as modular categories where $\Gamma$ is a finite $p$-group.

We refer to Appendix \ref{bcce} for an introduction to central braided crossed categories. A (faithfully graded) central braided $\Gamma$-crossed category $\D=\bigoplus\limits_{\g\in\Gamma}\C_\g$ with identity component $\M(A,\theta)$ (where $(A,\theta)$ is some metric group) is pointed by Remark \ref{centralpointed}. Hence, the set of isomorphism classes of simple objects in $\D$ is a group which is a central extension of $\Gamma$ by $A$, the group of isomorphism classes of simple objects of $\C_1=\M(A,\theta)$. The class of the central extension gives us an element $f\in H^2(\Gamma, A)$. Now, we can start from a central extension of $\Gamma$ by $A$ corresponding to an $f\in H^2(\Gamma, A)$ and ask whether there exist central braided $\Gamma$-crossed categories $\D$ giving rise to this central extension. According to \cite[\S 8]{ENO}, such central braided crossed categories exist when a certain obstruction $O_4(f)\in H^4(\Gamma,\Qlcl^*)$ vanishes, in which case they are classified up to equivalence by a certain $H^3(\Gamma,\Qlcl^*)$-torsor. We refer to \cite{ENO} for more on classification of braided crossed categories and in particular to {\it op. cit.} \S4 for the definition of various higher categorical groups (for example $\uuPic(\cdot), \uEqBr(\cdot)$) associated to braided fusion categories that play an important role in the classification. (See also Appendix \ref{bcce}.)

\brk
Since $\M(A,\theta)$ is a non-degenerate braided fusion category, we have an equivalence $\Theta:\uPic(\M(A,\theta))\cong\uEqBr(\M(A,\theta))$ (see \cite[\S5.4]{ENO} for details or Appendix \ref{bcce}). According to \cite[ \S8]{ENO}, $f\in H^2(\Gamma, A)$ classify morphisms of 1-groups $\Gamma\to \uPic(\M(A,\theta))\cong\uEqBr(\M(A,\theta))$ lifting the {\it trivial} group homomorphism $\Gamma\to \Pic(\M(A,\theta))\cong \EqBr(\M(A,\theta))$, or equivalently braided actions of $\Gamma$ on $\M(A,\theta)$ lifting the { trivial} action of $\Gamma$ on $(A,\theta)$. Such a morphism $\Gamma\to \uPic(\M(A,\theta))$ corresponding to $f\in H^2(\Gamma, A)$ can be lifted to a morphism of 2-groups $\Gamma\to \uuPic(\M(A,\theta))$ if and only if the obstruction $O_4(f)$ vanishes. In this case, the lifts are classified by a certain $H^3(\Gamma,\Qlcl^*)$-torsor. According to \cite{ENO} (see Theorem \ref{fusionandhomotopy}), the lifts $\Gamma\to \uuPic(\M(A,\theta))$ as above classify the central braided crossed categories corresponding to the $f\in H^2(\Gamma,A)$ up to equivalence.
\erk

Let $\tM$ be a central braided $\Gamma$-crossed category with identity component $\M_p^{anis}$ corresponding to an $f\in H^2(\Gamma, \f{F}_{p^2})$ (this implies that the corresponding obstruction $O_4(f)$ vanishes). Then we wish to construct a unipotent group $G$ with a Heisenberg idempotent $e$ such that $\tM_{G,e}\cong \tM$.

Let us now show that if we can realize one such central braided $\Gamma$-crossed category corresponding to $f$ from a unipotent group, then we can realize all. 

\bprop\label{torsor}
Let $\tM, \tM'$ be two central braided $\Gamma$-crossed categories corresponding to the class $f\in H^2(\Gamma, \f{F}_{p^2})$. Suppose there exists a unipotent group $G$ with a Heisenberg idempotent $e$ such that $\tM_{G,e}\cong \tM$. Then there exists a unipotent group $U$ with a Heisenberg idempotent $e'$  such that $\tM_{U,e'}\cong\tM'$.
\eprop

\bpf
Indeed $\tM'$ differs from $\tM$ in the associativity constraint by an element $\omega\in H^3(\Gamma,\Qlcl^{\times})$. By assumption, we have an extension $0\to H\to G\to \Gamma\to 0$ with a Heisenberg idempotent $e\in\D_G(G)$ corresponding to a Heisenberg admissible pair $(N,\L)$ for $G$ such that $\tM\cong \tM_{G,e}$. By \S\ref{cp+}, we also have a connected commutative unipotent group $V$ equipped with a $\Gamma$-action and a $\Gamma$-equivariant multiplicative local system $\K$ on $V$ and an extension $0\to V\to W\to \Gamma\to 0$ corresponding to the 3-cocycle $\omega$ (namely such that $\tM_{W,\K[2\dim V](\dim V)}\cong \Vec_\Gamma^{\omega}$). Let us consider the product of the group extensions: $$0\to H\times V\to U\to \Gamma \to 0.$$
Now we can check that $(N\times V, \L\boxtimes \K)$ is a Heisenberg admissible pair for $U$ and that the corresponding braided $\Gamma$-crossed category is $\tM'$.
\epf

Hence in order to complete the proof of Theorem \ref{conj1}(ii), given $f\in H^2(\Gamma, \F_{p^2})$, or equivalently a central extension $0\to \F_{p^2}\to K_{\Gamma}\to \Gamma\to 0$ (such that the obstruction $O_4(f)$ vanishes) it is enough to construct a unipotent group $U$ with $\pi_0(U)=\Gamma$ along with a Heisenberg idempotent such that the corresponding braided $\Gamma$-crossed category is one corresponding to the given central extension. (See \S\ref{completion} for details.)

\subsubsection{The central extension and obstruction in two step nilpotent case} \label{tce}
Let $H$ be a two step nilpotent unipotent group such that its center $N$ is connected and $H/N$ is commutative. Let $\Gamma$ be a finite $p$-group with an outer action on $H$ given by a group homomorphism $\Gamma\to \Out(H)$. This induces a braided action of $\Gamma$ on the braided monoidal category $\D_H(H)$. The outer action also induces an honest action of $\Gamma$ on the center $N$ as well as on $H/N=\Inn(H)$ since we have assumed that $\Inn(H)$ is commutative. Let us assume that the induced action of $\Gamma$ on $H/N$ is trivial. 

Let $\L$ be a $\Gamma$-equivariant multiplicative local system on $N$, such that $(N,\L)$ is a Heisenberg admissible pair on $H$. This means that the induced ($\Gamma$-equivariant) skew-symmetric biextension map $\phi_\L:(H/N)_{perf}\to (H/N)_{perf}^*$ is an isogeny whose kernel is denoted $K_\L$ as before. We have the corresponding non-degenerate quadratic form $\theta$ on $K_\L$ associated to this skew-symmetric isogeny. The admissible pair $(N,\L)$ gives us the Heisenberg idempotent $e=\L[2\dim N](\dim N)$ on $H$. 

As before, let $K\subset H$ be such that $K/N=K_\L$. For $h\in H$, let $c_h:H\to N$ be the ``commutator with $h$'' map: $H\ni g\mapsto ghg^{-1}h^{-1}\in N$. Recall that by Remark \ref{ssi} (see also \cite[Prop. 4.5]{De}), 
\beq\label{Kdef}
K=\{k\in H|c_k^*\L\cong \Qlcl\}. 
\eeq

Since $\L$ is $\Gamma$-equivariant, we get a braided action of  $\Gamma$ on $e\D_H(H)$ and hence on $\M_{H,e}\cong\M(K_\L,\theta)$. Since the action of $\Gamma$ on $H/N$ is trivial, we see that the braided action of $\Gamma$ on $\M_{H,e}$ induces the trivial action on $K_\L$.  According to \cite{ENO} this braided action 
corresponds to a certain central extension of $\Gamma$ by $K_\L$ as described in the previous section. Let us construct this braided action and central extension explicitly.

Let us pullback the exact sequence $0\to \Inn(H)\to \Aut(H)\to \Out(H)\to 0$ along the outer action map and obtain the central (since induced action of $\Gamma$ on $\Inn(H)$ is trivial) extension $$0\to \Inn(H)\to \Aut_{\Gamma}(H)\to \Gamma \to 0.$$ For $\gamma\in \Gamma$, let $\Aut_\gamma(H)\subset \Aut(H)$ be the $\Inn(H)$-coset that describes the outer action of $\gamma$ on $H$. In particular, $\Aut_1(H)=\Inn(H)=H/N$. For $\g\in\Gamma$ and $\tg\in \Aut_\g(H)$, let $c_{\tg}:H\to N$ be the map $c_{\tg}(h)=h\tg(h^{-1})\in N$. ($c_{\tg}(h)\in N$ since $\Gamma$ acts trivially on $H/N$.) In particular, for $Nh\in \Aut_1(H)=H/N$, $c_{Nh}:H\to N$ is the commutator with $h$ map $c_h$ defined previously. For each $\tg\in \Aut_\g(H)$, $c_{\tg}$ is a group homomorphism. Hence $c_{\tg}^*\L$ is a multiplicative local system on $H$ that is trivial when restricted to $N$ since $\L$ is $\Gamma$-equivariant. For $\g\in\Gamma$, let
\beq\label{kgamma}
K_\g=\{\tg\in \Aut_\g(H) | c_{\tg}^*\L\cong \Qlcl\}. 
\eeq
By (\ref{Kdef}), we see that $K_1=K/N=K_\L$.
 We have the relation $$c_{\tga\tgb}(h)=c_{\tga}(h)\ga(c_{\tgb}(h)).$$

Let $K_{\Gamma}$ be the disjoint union of all the $K_\g$ for $\g\in\Gamma$. Then we can check that in fact $K_\Gamma$ is a group and we have a central extension $$0\to K_1=K_\L\to K_{\Gamma}\to \Gamma\to 0.$$ $K_\Gamma$ has an honest action on $H$ and it induces the given outer action of $\Gamma$ on $H$. We also have a braided action of $K_\Gamma$ on $e\D_H(H)$ and $\M_{H,e}$ and the action of $K_\L\subset K_\Gamma$ corresponds to the conjugation action by the simple objects of $\M_{H,e}$. Hence this induces the braided action of the quotient $\Gamma$ on $\M_{H,e}$.

Thus from the outer action of $\Gamma$ on $H$ and the $\Gamma$-equivariant multiplicative local system $\L$ on the center $N$, we get a two cocycle $f:\Gamma\times \Gamma\to K_\L\subset H/N=\Inn(H)$ corresponding to the central extension constructed above by choosing lifts $\tg\in K_\g\subset \Aut_\g(H)\subset \Aut(H)$ for each $\g\in \Gamma$ such that $$\tga\tgb=f(\ga,\gb)\widetilde{\ga\gb}.$$ 

\noin Let $F:\Gamma\times \Gamma\to K\subset H$ be a lift of $f$ to $K\subset H$. Hence we get a 3-cocycle $\omega:\Gamma^3\to N$ such that we have 
\beq\label{omega}
\tilde{a}(F(b,c))F(a,bc)=\omega(a,b,c)F(a,b)F(ab,c), \mbox{ for } a,b,c\in \Gamma.
\eeq
\brk
Thus starting from an outer action of $\Gamma$ on $H$ we get an $\omega\in H^3(\Gamma,N)$. This is the obstruction to the existence of an extension of $\Gamma$ by $H$ corresponding to given outer action (see \cite[\S7]{EM}). 
\erk

\brk\label{centext}
Furthermore, if we have a $\Gamma$-equivariant Heisenberg admissible pair $(N,\L)$ as above, then we get a braided action of $\Gamma$ on $\M_{H,e}$ and the corresponding central extension $0\to K_\L\to K_\Gamma\to \Gamma\to 0$.
\erk

This data gives rise to an obstruction $\beta\in H^4(\Gamma,\Qlcl^*)$ to the existence of a central braided $\Gamma$-crossed category with trivial component $\M_{H,e}$. On the other hand we have the boundary map $\delta:H^3(\Gamma,N)\to H^4(\Gamma,\Qlcl^*)$ coming from the $\Gamma$-equivariant multiplicative local system $\L$ on $N$ (which can be thought of as a short exact sequence of $\Gamma$-modules $0\to \Qlcl^*\to \tN\to N\to 0$). Next we show that $\delta\omega= \beta$.

\subsubsection{The cohomological obstructions}
We continue to use the notations of the previous section in which we described a 3-cocycle $\omega\in H^3(\Gamma, N)$ which is the obstruction to the existence of an extension of $\Gamma$ by $H$ corresponding to the given outer action. We have a $\Gamma$-equivariant multiplicative local system $\L$ on $N$ which gives rise to a boundary map $\delta: H^3(\Gamma,N)\to H^4(\Gamma,\Qlcl^*)$.  We note that the outer action respected the Heisenberg admissible pair $(N,\L)$ on $H$ and that the induced action of $\Gamma$ on $H/N$ (and hence on $K_\L\subset H/N$) was trivial. The outer action defines a braided action of $\Gamma$ on $\M_{H,e}$ and we constructed the corresponding central extension $$0\to K_\L\to K_\Gamma \to \Gamma\to 0.$$  From this data, we get a $\beta\in H^4(\Gamma, \Qlcl^*)$ which is the obstruction to the existence a central braided $\Gamma$-crossed category corresponding to the braided action of  $\Gamma$ on $\M_{H,e}$. In this situation:
\bprop\label{bdry}
We have $\beta=\delta\omega$.
\eprop
\bpf
Let us first describe the boundary map $\delta$ explicitly. Let $\omega\in H^3(\Gamma,N)$ be represented by a 3-cocycle $\omega:\Gamma\times\Gamma\times\Gamma\to N$. For $a,b,c\in \Gamma$, let us choose a trivialization of the stalk $\L_{\omega(a,b,c)}$. We have $\omega(ab,c,d)^{-1}\cdot\omega(a,b,cd)^{-1}\cdot{}^{\ta}\omega(b,c,d)\cdot\omega(a,bc,d)\cdot\omega(a,b,c)=1$ for $a,b,c,d\in \Gamma$. Hence using our chosen trivializations of $\L_{\omega(\cdot,\cdot,\cdot)}$ and the $\Gamma$-equivariant multiplicative structure of the local system $\L$, we get an automorphism of the stalk $\L_1=\Qlcl$. This automorphism is $\delta\omega(a,b,c,d)$. 

We now describe the obstruction $\beta\in H^4(\Gamma,\Qlcl^*)$ explicitly. By Proposition \ref{cdo} (see also \cite[\S 4.3]{De}), the simple objects of $\M_{H,e}$ can be described as the translates $e^k$ of $e$ by $k\in K\subset H$ and moreover, for $k_1, k_2\in K$ we can identify $e^{k_1}\ast e^{k_2}$ with $e^{k_1k_2}$. The isomorphism class of $e^k$ only depends on the coset $Nk$. For $n\in N$ and any $k\in K$, a choice of trivialization of the stalk $\L_n$ gives us an isomorphism $e^{k}\to e^{nk}$. For $a,b,c\in \Gamma$, we have a chosen trivialization of the stalk $\L_{\omega(a,b,c)}$ as in the previous paragraph. 

Note that for $a\in \Gamma$, we have chosen lifts $\ta\in K_a\subset \Aut_a(H)$ (see \S\ref{tce}). We define the corresponding braided automorphism (also denoted by $\ta$) of $\M_{H,e}$ as pullback via $\ta^{-1}:H\to H$.  Let $\M_{H,e,\ta}$ denote the corresponding braided $\M_{H,e}$-bitorsor which is equal to $\M_{H,e}$ as a left-torsor and where the right action is twisted by the braided automorphism $\ta$. For $X\in\M_{H,e}$ and $a,b\in \Gamma$, $\ta\tb(X)=e^{F(a,b)}\ast\tab(X)\ast e^{F(a,b)^{-1}}$. We will define a system of products (see \cite[\S8.6]{ENO}) $\M_{H,e,\ta}\boxtimes_{\M_{H,e}} \M_{H,e,\tb}\to \M_{H,e,\tab}$ mapping $(X,Y)\mapsto X\ast\ta(Y)\ast e^{F(a,b)}$ where we have used the identifications of all bitorsors with $\M_{H,e}$ as left-torsors. Now for $a,b,c\in \Gamma$, we have two functors $\M_{H,e,\ta}\times \M_{H,e,\tb}\times \M_{H,e,\tc}\to \M_{H,e,\t{abc}}$ corresponding to the two bracketings. They are given by $$(X,Y,Z)\mapsto X\ast\ta(Y)\ast\ta\tb(Z)\ast e^{{ }^{\ta} F(b,c)F(a,bc)},$$ $$(X,Y,Z)\mapsto X\ast\ta(Y)\ast\ta\tb(Z)\ast e^{F(a,b)F(ab,c)}.$$ Hence by (\ref{omega}), the trivialization of the stalk $\L_{\omega(a,b,c)}$ induces a natural isomorphism $\zeta_{a,b,c}$ between these two functors. Now for $a,b,c,d\in \Gamma$ we have five functors $\M_{H,e,\ta}\times \M_{H,e,\tb}\times \M_{H,e,\tc}\times \M_{H,e,\td}\to \M_{H,e,\t{abcd}}$ corresponding to the five bracketings and a pentagon formed by the natural isomorphisms between these functors defined using the $\zeta_{\cdot,\cdot,\cdot}$ defined above. Let $\beta(a,b,c,d)\in\Qlcl^*$ denote the natural automorphism of one of the five functors obtained from the above pentagon (namely $\beta$ measures the failure of the commutativity of the pentagon). It follows that $\beta$ is a 4-cocycle and that the associativity constraint $\zeta$ can be modified to satisfy the pentagon axiom if and only if $\beta$ is cohomologically trivial, or in other words, $\beta$ is the obstruction to the existence of the central braided crossed category. Comparing with the explicit computation of $\delta\omega$ above, we see that $\beta=\delta\omega$.

\epf

\subsubsection{Constructing an outer action from a central extension}
Let $\Gamma$ be a finite $p$-group and let $$0\to \f{F}_{p^2}\to K_{\Gamma} \to \Gamma \to 0$$ be a central extension. We will now reverse the construction of \S\ref{tce}.
Namely, we will construct a two step nilpotent group $H$ with connected center $N$ and an embedding $\F_{p^2}\subset H/N$ with an action of $K_\Gamma$ such 
that $\F_{p^2}\subset H/N = \Inn(H)$ acts by its natural 
conjugation action on $H$ (this induces an outer action of $\Gamma$ on $H$) and a $\Gamma$-equivariant 
Heisenberg admissible pair $(N,\L)$ as described in \S\ref{tce}, such that the corresponding metric group is the anisotropic 
metric group $(\F_{p^2},\theta=\zeta^{Nm(\cdot)})$ and such that the corresponding extension of $\Gamma$ by $K_\L=\f{F}_{p^2}$ obtained using the construction of \S\ref{tce} coincides with $K_\Gamma$.

Let us fix a positive integer $n$ and let $\W$ be the ring scheme of Witt vectors of length $n$. Let 
$\tN=Maps(\Gamma, \W)$, or the group algebra of $\Gamma$ over the ring $\W$. Then $\tN$ is a connected commutative 
unipotent group with a left $\Gamma$-action by right translations of functions. Note that we have 
$H^i(\Gamma, \tN)=0$ for $i>0$ by Shapiro's lemma.

We have the embedding 
$C:=\f{Z}/{p^n}\f{Z}\subset \W{\hookrightarrow}\tN$ of constant maps 
which is a $\Gamma$-equivariant map with $\Gamma$ acting trivially on $\W$. 
Let us form the quotient $N:=\tN/C$ which is a connected commutative unipotent group equipped with the induced action of $\Gamma$. 

\blem\label{embed}
We have the following structure:
\bit
\item[(i)] A $\Gamma$-equivariant multiplicative local system $\L$ on $N$.
\item[(ii)] Embeddings $\G_a\hookrightarrow\W\subset N$ such that the induced $\Gamma$-action is trivial on $\W$ and $\L$ restricts to the Artin-Schreier local systems on $\W$ and $\G_a.$
\eit
\elem
\bpf
Consider the $\Gamma$-equivariant central extension $$0\to C\to\tN\to N\to 0.$$ 
This gives us (i), the $\Gamma$-equivariant multiplicative local system $\L$ on $N$ 
(after choosing an identification $C\cong\mu_{p^n}\subset \Qlcl^{*}$). Restricting the above central extension to the constants $\W\subset \tN$,
we get the central extension $0\to C\to \W\to\W/C\to 0$, which can be identified with the Artin-Schreier central extension $0\to C\to \W{\xrightarrow{w\mapsto F(w)-w}}\W\to 0$, where $F$ is the relative Frobenius map. This gives us an identification $\W\cong\W/C$ and an embedding $\W\cong\W/C\subset \tN/C=N$ and we see that under this embedding $\L|_{\W}$ is the Artin-Schreier multiplicative local system. The action of
$\Gamma$ on $\W$ is trivial under this embedding. Also, we have $\G_a\hookrightarrow \W$ 
under which the Artin-Schreier local system on $\W$ restricts to the Artin-Schreier local system on $\G_a$. 
\epf

An immediate consequence is the following:
\blem\label{defB}
Let $B:\G_a\times \G_a\to \G_a\subset N$ be the skew-symmetric biadditive 
map given by $B(\l_1,\l_2)=\l_1{\l_2}^p-{\l_1}^p\l_2\in \G_a\subset N$. The pullback $B^*\L$  induces a 
skew-symmetric biextension $\Gap\to (\Gap)^*$ which is an isogeny with kernel $\F_{p^2}$ and the induced 
metric is the anisotropic metric.
\elem
\bpf
This follows from Example \ref{ex} since $\L|_{\G_a}$ is the Artin-Schreier multiplicative local system and $B$ identifies with the commutator map of Example \ref{ex}. 
\epf

For $\lambda\in \G_a$, 
let $\tN_\lambda:=\{\phi\in Maps(K_{\Gamma}, \W)| \phi(kx)=\phi(x)+k^p\lambda\ \hbox{ for each }k\in\F_{p^2}, x\in K_\Gamma\}$, the space of $\F_{p^2}$-quasi-invariant maps. Note that
$\tN_0$ is the set of $\F_{p^2}$-invariant maps and hence can be identified with $\tN=Maps(\Gamma,\W)$. We have the natural action of $K_{\Gamma}$
on $\tN_\lambda$ given by $y\phi(x)=\phi(xy)$. For 
$k\in \F_{p^2}\subset K_\Gamma$, its action on $\tN_\lambda$ is given by $k\phi(x)=\phi(xk=kx)=\phi(x)+k^p\lambda.$
Each $\tN_\lambda$ is a $\tN$-bitorsor and addition induces 
$\Gamma$-equivariant maps $\tN_{\lambda_1}\times \tN_{\lambda_2}\to \tN_{\lambda_1+\lambda_2}$ for $\lambda_i\in \G_a.$

For $\lambda\in \G_a$, define $N_\lambda:=\tN_{\lambda}/C$. In particular $N_0=N$. 
For $\phi\in \tN_\lambda,$ let $\bar{\phi}$ denote its image in $N_\lambda.$ Note that for a constant function
$w\in\W\subset \tN, \bar{w}=F(w)-w$  considered as an element of the $\W$ embedded in $N$ (see Lemma \ref{embed}), where $F:\W\to \W$ 
is the relative Frobenius.
We have the induced action of $K_\Gamma$ on each $N_\lambda$ and the induced $\Gamma$-equivariant addition
maps $N_{\lambda_1}\times N_{\lambda_2}\to N_{\lambda_1+\lambda_2}$ for $\lambda_i\in \G_a$. For 
$k\in \F_{p^2}\subset K_\Gamma$, its action on $N_\lambda$ is given by 
$$k\bar{\phi}=\bar{\phi+k^p\lambda}=\bar{\phi}+F(k^p\lambda)-k^p\lambda$$ 
$$=\bar{\phi}+k^{p^2}\lambda^p-k^p\lambda=\bar{\phi}+k\lambda^p-k^p\lambda$$
$$=\bar{\phi}+B(k,\lambda).$$

Let us now define the unipotent group $H$. The underlying variety $H=\{(\lambda, f)|\lambda\in\G_a, f\in N_\lambda\}$.
Define a product on $H$ by setting 
$$(\lambda_1,f_1)\cdot(\lambda_2,f_2)=(\lambda_1+\lambda_2, f_1+f_2+\l_1{\l_2}^p)).$$
(Here $\l_1{\l_2}^p\in \G_a\subset N$.) Note that $(\lambda, f)^{-1}=(-\l, \l^{p+1}-f)$.  The center of $H$ is $N$ and $\Inn(H)=H/N=\G_a$ with the inner action of $\mu\in\G_a$ given by ${}^{\mu}(\l,f)=(\l,f+B(\mu,\l))$ where $B$ is as defined in Lemma \ref{defB}. The commutator map $[\cdot,\cdot]:H\times H\to N\subset H$ is given by $[(\l_1,f_1),(\l_2,f_2)]=B(\l_1,\l_2).$ By this, Remark \ref{ssi} and Lemma \ref{defB} we have:

\blem\label{fp2}
The pair $(N,\L)$ is a Heisenberg admissible pair on $H$ and the corresponding metric group is the anisotropic metric group $(\F_{p^2},\theta)$.
\elem

The action of $K_\Gamma$ on each $N_\l$ induces an action of $K_\Gamma$ on the group $H$ defined by $x(\l,f)=(\l,xf)$. The induced action of $\Gamma$ on $N$ is the original action. The induced action on $H/N$ is trivial. As before, the action of $k\in\F_{p^2}\subset K_\Gamma$ is given by $k(\l,f)=(\l,f+B(k,\l)).$ We see that this identifies with the inner action of $k\in\F_{p^2}\subset \G_a=H/N.$  

We can now prove:
\bprop
The central extension obtained from the data as above (using the construction of \S\ref{tce}) is the central extension $0\to \f{F}_{p^2}\to K_{\Gamma} \to \Gamma \to 0$ that we started out with.
\eprop
\bpf
It is clear (for example by Lemma \ref{fp2}) that in this case $K_\L=\F_{p^2}$. For $y\in K_\Gamma$, let $c_y:H\to N$ be defined by $h=(\l,\bar{\phi})\mapsto h\cdot y(h^{-1})=\bar{\phi}-y\bar{\phi}.$ This is a group homomorphism which can in fact be lifted to a group homomorphism $C_y:H\to \tN$ by $(\l,\bar{\phi})\mapsto \phi-y\phi.$ One can readily verify that this is indeed a well defined homomorphism. This means that $c_y^*\L$ is the trivial multiplicative local system on $H$. Hence using (\ref{kgamma}) we see that $K_\Gamma$ is indeed the central extension obtained by the construction of \S\ref{tce}.
\epf

\subsubsection{Completion of the proof}\label{completion}
Let $\C\in\mathfrak{C}_p^-$. By Proposition \ref{dgno2},  $\C\cong \tM'^\Gamma$ as modular categories, where $\tM'$ is  a (positive spherical) central braided $\Gamma$-crossed category with trivial component $\M_p^{anis}$ and $\Gamma$ is a finite $p$-group. The group of isomorphism classes of simple objects in the central braided crossed category $\tM'$ gives us an $f\in H^2(\Gamma, \F_{p^2})$ such that the corresponding obstruction $\beta\in H^4(\Gamma, \Qlcl^*)$ vanishes. Note that $\beta$ can be defined as an element of $H^4(\Gamma,\mu_p)$. Hence $\beta$ must vanish in some $H^4(\Gamma, \mu_{p^n})$. Let us fix one such $n$. Let $H$ be the two-step nilpotent unipotent group  equipped with the outer action of $\Gamma$ and the $\Gamma$-equivariant admissible pair $(N,\L)$ corresponding to $f\in H^2(\Gamma,\F_{p^2})$ constructed in the previous section (for the chosen integer $n$). Let $\omega\in H^3(\Gamma,N)$ be the obstruction to the existence of a group extension corresponding to the outer action of $\Gamma$ on $H$ constructed above. The multiplicative local system $\L$ was defined using the $\Gamma$-equivariant central extension $0\to \mu_{p^n}\cong C\to \tN\to N\to 0$. Since $H^i(\Gamma,\tN)=0$ for positive integers $i$, the boundary map defines an isomorphism $\delta: H^3(\Gamma,N)\cong H^4(\Gamma,\mu_{p^n})$. Since $\beta=\delta\omega=0$, we conclude that $\omega$ vanishes and hence there exists an extension $0\to H\to G\to \Gamma\to 0$. It follows that $(N,\L)$ is a special Heisenberg admissible pair for $G$ and that the corresponding central braided $\Gamma$-crossed category $\tM_{G,e}$ has identity component $\M_p^{anis}$ and group of isomorphism classes of simple objects is the central extension corresponding to $f$. 

Hence $\tM'$ and $\tM_{G,e}$ both give rise to the same central extension $f$. By Proposition \ref{torsor}, there exists a unipotent group $U$ with a Heisenberg idempotent $e'$ such that $\tM'\cong \tM_{U,e'}$ as spherical central braided crossed categories. Hence $\M_{U,e'}\cong\tM_{U,e'}^\Gamma\cong \C$ as modular categories.

By \S\ref{dtoc}, we can pass to a connected unipotent group. Hence the proof of Theorem \ref{conj1} is now complete.

\appendix

\section{Appendix: Skew-symmetric biextensions}\label{A}

In this appendix, we prove some general results about connected commutative (perfect) unipotent groups and skew-symmetric biextensions. (See \cite[Appendix A.6-11]{B},\cite[\S2]{Da} for more on skew-symmetric biextensions.) Since it is convenient to talk about skew-symmetric biextensions in the setting of perfect schemes, we will work in the category of perfect unipotent groups over $\k$, even though many of the results and arguments also work in the category of unipotent groups. Let $\cpu$ denote the category of commutative perfect unipotent groups over $\k$ and let $\cpuc$ be the full subcategory of connected commutative perfect unipotent groups. The category $\cpu$ is abelian and $\cpuc$ is an exact subcategory. Serre duality defines an involution $\cpuc\ni V\mapsto V^*\in \cpuc$ which is an exact anti-autoequivalence of $\cpuc$. For $V\in \cpuc$, we think of $V^*\in \cpuc$ as the space parametrizing the multiplicative local systems on $V$ (after fixing an embedding $\Q_p/\f{Z}_p\hookrightarrow \Qlcl^*$). We refer to \cite[Appendix A]{B} for more on the notions of Serre duality and biextensions.

\subsection{Group actions on a connected commutative unipotent group}
Let $\Gamma$ be a group. Let $V\in\cpu$ (resp. $\in \cpuc$) be equipped with an action of $\Gamma$. Let $\Gcpu$ (resp. $\Gcpuc$) denote the category whose objects are $V$ (resp. connected $V$) equipped with a $\Gamma$-action and whose morphisms are $\Gamma$-equivariant maps of commutative perfect unipotent groups. Then $\Gcpu$ is an abelian category and $\Gcpuc$ is an exact subcategory. We prove some results in this situation which reduce to basic well-known results in the case when $V$ is a vector space over $\k$ and $\Gamma$ acts by linear transformations.

\bdefn
An object $V\in \Gcpuc$ is called {\it simple} if it is nonzero and has no proper {\it connected} $\Gamma$-invariant subgroups.
\edefn

It is easy to see that:
\blem
Every $V\in\Gcpuc$ has a filtration $\{0\}=V_0\subset V_1\subset \cdots\subset V_{n-1}\subset V_n=V$ by connected $\Gamma$-invariant subgroups such that each successive quotient $V_{i+1}/V_i$ is simple.
\elem

We see that a morphism in $\Gamma$-$\cpuc$ to a simple object is either zero or surjective and that a morphism from a simple object is either zero or has finite kernel. Hence we have the following:

\blem
Let $V,W$ be simple objects in $\Gcpuc$ and $f:V\to W$ a morphism in $\Gcpuc$. Then $f$ is either zero or an isogeny.
\elem

\blem
Let $V\in \Gcpuc$ be simple. Then $p\cdot V=0$, i.e. as a connected commutative perfect unipotent group, $V\cong \Gap^r$ for some positive integer $r$.
\elem
\bpf
Consider the morphism  $p\cdot\Id_V: V\to V$ in $\Gcpuc.$ Since some power of $p$ annihilates $V$, the map $p\cdot\Id_V$ is not an isogeny. Hence it must be zero.
\epf

\blem
Let $V\in\Gcpuc$ be simple. Let $\g\in \Gamma$ be a central element whose order is a power of $p$. Then $\g$ acts trivially on $V$.
\elem
\bpf
Suppose $\g^{p^n}=1$. Since $\g\in \Gamma$ is central, $\g:V\to V$ is a morphism in $\Gcpuc$ and so is $(\g-\Id)$. Since V is simple $p\cdot V=0$. Hence $(\g-\Id)^{p^n}=\g^{p^n}-\Id=0.$ Hence $\g-\Id$ cannot be an isogeny. Hence $\g=\Id$ since $V$ is simple.
\epf

\brk
If $\Gamma$ is a (possibly infinite) $p$-group and $V\in\Gcpuc$ is simple, then the center $Z(\Gamma)$ acts trivially on $V$.
\erk

\bprop
If $\Gamma$ is a (possibly infinite) nilpotent $p$-group, then there is a unique simple object in $\Gcpuc$, namely $\Gap$ equipped with the trivial action of $\Gamma$.
\eprop
\bpf
We prove this by induction on the length of the ascending central series of $\Gamma$. By the remark, the proposition holds for commutative $p$-groups $\Gamma$. ($Z(\Gamma)=\Gamma$ must act trivially on a simple $V$, and consequently $V$ must be one dimensional.) Now suppose $\Gamma$ is such that the proposition holds true for all $p$-groups of lower nilpotence class. Let $V\in\Gcpuc$ be simple. By the remark above, $Z(\Gamma)$ acts trivially on $V$.
Hence the $\Gamma$-action on $V$ comes from an action of the quotient $\Gamma/Z(\Gamma)$. $V$ is simple for this action. Hence by our inductive hypothesis, $V\cong \Gap$ with trivial action of $\Gamma/Z(\Gamma)$. Hence the action of $\Gamma$ is trivial as well.
\epf

\bcor\label{keylemma}
Let $\Gamma$ be a nilpotent $p$-group. Let $V\in\Gcpuc$ with $\dim(V)=n$. Then there is a filtration $\{0\}=V_0\subset V_1\subset \cdots\subset V_{n-1}\subset V_n=V$ by connected $\Gamma$-invariant subgroups such that each $V_{i+1}/V_i\cong \G_a$ with $\Gamma$ acting trivially.
\ecor

\subsection{Invariant isotropic subgroups}
As before, let $\Gamma$ be a group and let $V\in\Gcpuc$. Let $\phi:V\to V^*$ be a $\Gamma$-equivariant skew-symmetric biextension. (See \cite[Appendix A.10]{B}.) We can think of such a skew-symmetric biextension as a bimultiplicative local system on $V\times V$ whose restriction to the diagonal is trivial. Let $W\subset V$ be a connected subgroup. Let $W^\perp (\in \cpu)$ be the kernel of the composition $V\stackrel{\phi}{\to} V^*\onto W^*$. We say that $W$ is isotropic if $W\subset W^{\perp}$. In this section, we prove that if $\dim(V)\geq 2$ and if $\Gamma$ is a nilpotent $p$-group, then there exists a non-trivial $\Gamma$-invariant isotropic subgroup.

\blem
Let $n$ be the smallest integer such that $p^n\cdot V=0$. Let $r$ be an integer such that $\frac{n}{2}\leq r \leq n-1.$ Then $p^r\cdot V$ is a proper $\Gamma$-invariant isotropic subgroup.
\elem
\bpf
It is clear that $p^r\cdot V$ is a proper $\Gamma$-invariant subgroup for each positive integer $r\leq n-1$. For $r\geq \frac{n}{2}$, $p^r\cdot(p^r\cdot V)=0$, hence $p^r\cdot (p^r\cdot V)^*=0$. Hence we see that the composition $$p^r\cdot V\subset V\to V^*\onto (p^r\cdot V)^*$$ is zero. ($p^r\cdot x\mapsto p^r\cdot \phi(x)\mapsto p^r\cdot \phi(x)|_{p^r\cdot V}=0$ in $(p^r\cdot V)^*$ which is annihilated by $p^r$.)
\epf

\bprop
Let $\Gamma$ be a nilpotent $p$-group, $V\in\cpuc$ and let $\phi:V\to V^*$  be a $\Gamma$-equivariant skew-symmetric biextension. If $\dim(V)\geq 2$, then there exists a nonzero $\Gamma$-invariant isotropic subgroup in $V$.
\eprop
\bpf
If $p\cdot V\neq 0$, then the previous lemma gives us a proper $\Gamma$-invariant isotropic subgroup in $V$. Hence let us assume $p\cdot V=0$, in which case the previous lemma does not tell us anything. Let us use Corollary \ref{keylemma}. Since $\dim(V)\geq 2$, there exists a two-dimensional $\Gamma$-invariant subgroup $V_2\subset V$. Since $p\cdot V_2 =0, V_2\cong \Gap^2.$ Consider the induced $\Gamma$-equivariant skew-symmetric biextension of 
$V_2$. Let $V_1\subset V_2$ be a one-dimensional $\Gamma$-fixed subgroup guaranteed by Corollary \ref{keylemma}. If $V_1$ is isotropic in $V_2$, then it is also isotropic in $V$ and we are done. If not, then its perpendicular $V_1^{\perp_{V_2}}$ in $V_2$ must be one-dimensional and $V_1+(V_1^{\perp_{V_2}})^\circ=V_2.$ Since $V_1\subset V_2$ is $\Gamma$-stable, $(V_1^{\perp_{V_2}})^\circ$ must also be $\Gamma$-invariant and since it is one dimensional, the action of $\Gamma$ on $(V_1^{\perp_{V_2}})^\circ$ must be trivial and consequently action of $\Gamma$ on $V_2$ must be trivial. Now by \cite[Lemma A.29]{B}, $V_2$ has a one-dimensional isotropic subgroup $W$. Since the action of $\Gamma$ on $V_2$ is trivial, $W$ is $\Gamma$-stable.
\epf

\bcor\label{mgiis}
Let $\Gamma$ be a nilpotent $p$-group and let $V\in \Gcpuc$ be equipped with a $\Gamma$-equivariant skew-symmetric biextension $V\stackrel{\phi}{\to}V^*$. Then there exists a $\Gamma$-invariant maximal isotropic subgroup $W\subset V$ and for any such $W$, $\dim(W^{\perp})\leq \dim(W)+1$ and $W$ is also maximal among all (not necessarily $\Gamma$-invariant) isotropic subgroups of $V$.
\ecor
\bpf
Clearly maximal $\Gamma$-invariant isotropic subgroups exist. Let  $W$ be any $\Gamma$-invariant isotropic subgroup of $V$. Then $(W^\perp)^\circ$ is also $\Gamma$-stable and we have the induced $\Gamma$-equivariant skew-symmetric biextension of $(W^{\perp})^\circ/W$. $\Gamma$-invariant isotropic subgroups of $V$ containing $W$ are in a natural bijective correspondence with $\Gamma$-invariant isotropic subgroups of $(W^{\perp})^\circ/W$. If $W$ is a maximal $\Gamma$-invariant isotropic subgroup, then $(W^{\perp})^\circ/W$ has no nonzero $\Gamma$-invariant isotropic subgroups. Hence this space has dimension less than or equal to 1.
\epf

\subsection{Metric groups associated to skew-symmetric biextensions}\label{mgassssbe}
\bdefn\label{defmg}
Let $A$ be a finite abelian group and let $\theta: A\to \Qlcl^*$ be a function. We say that $\theta$ is a quadratic form if $\theta(a)=\theta(-a)$ and if the symmetric function $b(a,b):=\frac{\theta(a+b)}{\theta(a)\theta(b)}$ is a bicharacter $b:A\times A\to \Qlcl^*$. We say that $\theta$ is non-degenerate if the bicharacter $b$ is non-degenerate. A pre-metric group is a pair $(A,\theta)$ where $A$ is a finite abelian group and $\theta:A\to \Qlcl^*$ is a quadratic form. If the quadratic form is non-degenerate, we say that $(A,\theta)$ is a metric group.
\edefn

Given a metric group $(A,\theta)$, we have the associated non-degenerate bicharacter $b$. A subgroup $B\subset A$ is said to be isotropic if $\theta|_B=1$. In this case, $B\subset B^\perp$. For an isotropic subgroup $B$, the subquotient $(B^\perp/B, \theta')$ is also a metric group, where $\theta'$ is the quadratic form on $B^\perp/B$ induced by $\theta$. A metric group is said to be {\it anisotropic} if it has no non-trivial isotropic subgroups. If $L$ is a maximal isotropic subgroup of a metric group $(A,\theta)$, then the subquotient $(L^\perp/L,\theta')$ is an anisotropic metric group whose isomorphism class is independent of the choice of maximal isotropic subgroup $L$. We say that the metric group $(A,\theta)$ is {\it Witt-equivalent} to the anisotropic metric group $(L^\perp/L,\theta')$.

Let us now give two examples of {\it anisotropic} metric $p$-groups which are of particular interest in the theory of character sheaves on unipotent groups (see Theorem \ref{swarny} below):
\bit
\item The trivial metric group $(0, 1)$.
\item $(\F_{p^2}, \zeta^{{Nm}(\cdot)})$ where $Nm: \F_{p^2}\to \Fp$ is the norm map and $\zeta\in \Qlcl^*$ is a primitive $p$-th root of unity. 
\eit

Now let $V$ be a connected commutative perfect unipotent group and let $\phi:V\to V^*$ be a skew-symmetric isogeny. According to \cite[Appendix A.10]{B}, \cite[\S2.3]{Da}, to such a skew-symmetric isogeny, we can associate a metric group $(A,\theta)$. Here $A$ is the kernel of $\phi$. Then we can define the associated non-degenerate quadratic form $\theta:A\to \Qlcl^*$ (see \cite[Appendix A.10]{B}, \cite[\S 2.3]{Da}). The following theorem is proved in \cite{Da}:
\bthm{(\cite{Da})}\label{swarny}
The metric group $(A,\theta)$ associated to the skew-symmetric isogeny $\phi:V\to V^*$ is Witt-equivalent to the trivial anisotropic metric group if $\dim(V)$ is even and is Witt-equivalent to the anisotropic metric group $(\F_{p^2}, \zeta^{Nm(\cdot)})$ if $\dim(V)$ is odd.
\ethm

\section{Appendix: The modular category associated to a metric group and braided crossed categories}\label{mcassmg}

\subsection{The modular category associated to a metric group}
Given a pre-metric group $(A,\theta)$ (see Definition \ref{defmg}), there exists a corresponding pointed braided category $\M(A,\theta)$ whose group of isomorphism classes of simple objects is isomorphic to $A$ and the braiding on the square of a simple object is given by $\theta$. This determines the braided monoidal category up to equivalence. We equip this braided pointed category with the positive spherical structure. With this structure, the twist (an automorphism of the identity functor on the pointed braided category) is given by $\theta$. We denote this pre-modular category by $\M(A,\theta)$. (See \cite[\S 2.11.4-5 and Appendix D]{DGNO} for more details.) If $(A,\theta)$ is a metric group, $\M(A,\theta)$ is a modular category.

The modular category defined by the metric group $(\F_{p^2},\zeta^{Nm(\cdot)})$ is of particular interest and we denote it by $\M_p^{anis}$.

\subsection{Braided crossed categories and equivariantization}\label{bcce}
Let $(A,\theta)$ be a metric group and let $\Gamma$ be a finite group. Let $\D$ be a {\it faithfully graded} braided $\Gamma$-crossed category with identity component $\M(A,\theta)$. This means that we have the following structure on $\D$:
\bit
\item A grading $\D=\bigoplus\limits_{\g\in\Gamma}\C_\g$ where each $\C_\g$ is non-zero and $\C_1\cong \M(A,\theta)$. 
\item A monoidal action of $\Gamma$ on $\D$ such that $\ga(\C_\gb)\subset \C_{\ga\gb\ga^{-1}}$ for each $\ga,\gb\in\Gamma$.
\item For $\g\in \Gamma, X\in\C_\g$ and $Y\in \D$ isomorphisms 
$$c_{X,Y}:X\otimes Y\stackrel{\cong}{\rightarrow} \g(Y)\otimes X$$ functorial in $X,Y$ and satisfying certain compatibility conditions.
\eit
For a braided $\Gamma$-crossed category $\D$, the equivariantization $\D^\Gamma$ has the structure of a braided monoidal category. We refer to \cite[\S 4.4.3]{DGNO} for a precise definition and properties of braided crossed categories and related concepts.

In particular, the monoidal action of $\Gamma$ on $\D$ restricts to an action of $\Gamma$ on the braided monoidal category $\M(A,\theta)$ by braided autoequivalences which in turn determines an action of $\Gamma$ on the metric group $(A,\theta)$.

\bdefn\label{central}
 We say that such a braided $\Gamma$-crossed category as above is {\it central} if the induced action of $\Gamma$ on the metric group $(A,\theta)$ is trivial.
\edefn

Let $\D=\bigoplus\limits_{\g\in\Gamma}\C_\g$ be a braided $\Gamma$-crossed category. Each $\C_\g$ is a $\C_1$-module category. In fact, by \cite[Thm. 6.1]{ENO}, it is an invertible (under tensor product of module categories, see {\it op. cit.} \S 3) $\C_1$-module category. For a braided fusion category $\C_1$, let $\uuPic(\C_1)$ denote the categorical 2-group whose objects are invertible $\C_1$-module categories under tensor product, 1-morphisms are equivalences of $\C_1$-module categories and 2-morphisms are isomorphisms of such equivalences. As proved in {\it op. cit.}, the assignment $\gamma\mapsto \C_\g$ along with the braided $\Gamma$-crossed structure of $\D$ defines a morphism of 2-groups $\Gamma\to \uuPic(\C_1)$. In fact, the following result is proved in {\it op. cit. \S7.8} 
\bthm\label{fusionandhomotopy}(\cite[Thm. 7.12]{ENO})
Equivalence classes of faithfully graded braided $\Gamma$-crossed categories with trivial component $\C_1$ are in bijection with morphisms of categorical 2-groups $\Gamma\to \uuPic(\C_1)$.
\ethm

We can truncate (by forgetting the highest level morphisms and identifying isomorphic lower ones) the 2-group $\uuPic(\C_1)$ to a 1-group $\uPic(\C_1)$ and further to an ordinary group $\Pic(\C_1)$ of equivalence classes of invertible $\C_1$-module categories. Let $\uEqBr(\C_1)$ denote the 1-group whose objects are braided autoequivalences of $\C_1$ and whose morphisms are isomorphisms of such braided autoequivalences. We can truncate this to get an ordinary group $\EqBr(\C_1)$. According to \cite[\S5.4]{ENO}, there is a monoidal functor $\Theta:\uPic(\C_1)\to \uEqBr(\C_1)$ which is an equivalence in case the braided fusion category $\C_1$ is non-degenerate. 

\brk
A morphism of 1-groups $\Gamma\to \uEqBr(\C_1)$ is equivalent to a braided action of $\Gamma$ on $\C_1$. Note that a morphism $\Gamma\to \uuPic(\C_1)$ gives us a morphism $\Gamma\to \uEqBr(\C_1)$ using the monoidal functor $\Theta$. In other words we get a braided action of $\Gamma$ on $\C_1$. This is precisely the braided action coming from the braided crossed category structure defined by the morphism $\Gamma\to \uuPic(\C_1)$.
\erk

Now let us consider a faithfully graded central braided $\Gamma$-crossed category $\D$ with trivial component $\M(A,\theta)$ (where $(A,\theta)$ is a metric group), corresponding to a morphism $\Gamma\to \uuPic(\M(A,\theta))$. Centrality means that the induced map $\Gamma\to O(A,\theta)\cong \EqBr(\M(A,\theta))\cong\Pic(\M(A,\theta))$ is trivial. Here $O(A,\theta)$ is the group of automorphisms of $A$ preserving the quadratic form $\theta$. Hence in this case, all the graded components $\C_\g$ are equivalent to the trivial module category $\M(A,\theta)$ as $\M(A,\theta)$-module categories. 
\brk\label{centralpointed}
In particular, this implies that $\D$ is pointed. 
\erk
Let us equip $\D$ with the positive spherical structure (see \cite[Ex. 2.26]{DGNO}). Since the identity component $\M(A,\theta)$ is a non-degenerate braided category and the grading on $\D$ is faithful, the equivariantization $\D^\Gamma$ is also a non-degenerate braided category (by \cite[Proposition 4.56(ii)]{DGNO}). Since $\D$ is equipped with a (positive) spherical structure, $\D^{\Gamma}$ also has the structure of a modular category.

\bprop\label{core}
Let $\D$ be faithfully graded braided $\Gamma$-crossed category with identity component $\M(A,\theta)$, where 
$\Gamma$ is a finite group and $(A,\theta)$ is a metric group. Let $\C$ be the braided monoidal category 
$\D^\Gamma$. Let $L\subset A$ be a $\Gamma$-invariant isotropic subgroup and form the subquotient metric group $(L^{\perp}/L,\theta')$. Then there exists a (faithfully graded) braided crossed category whose identity component is equivalent to $\M(L^{\perp}/L,\theta')$ and whose equivariantization is equivalent to the braided monoidal category $\C$.
\eprop
\bpf
The full subcategory of $\M(A,\theta)$ generated by the simple objects corresponding to $L\subset A$ with the restriction of associativity and the braiding is the symmetric monoidal category $\M(L,1=\theta|_L)\cong \Vec_L$ and is equipped with a fiber functor $\Vec_L\to \Vec$. The braided action of $\Gamma$ on $\M(A,\theta)$ induces one on $\Vec_L$ and after equivariantization, we get a symmetric monoidal category $\mathcal{E}:=\Vec_L^{\Gamma}\subset \D^\Gamma= \C$ equipped with the fiber functor $\E=\Vec_L^{\Gamma}\to \Vec_L\to \Vec$. The centralizer of $\E$ in $\C$ is the subcategory $\E'=\M(L^{\perp},\theta|_{L^{\perp}})^{\Gamma}$ of $\M(A,\theta)^\Gamma$. The de-equivariantization (see \cite[\S4.4.7,8]{DGNO}) $\C\boxtimes_{\E}\Vec$ is a faithfully graded braided crossed category with identity component $\E'\boxtimes_{\E}\Vec\cong\M(L^{\perp}/L,\theta)$ and whose equivariantization is equivalent to $\C$. (See \cite[Prop. 4.56]{DGNO}.) 
\epf

\brk\label{pgroup}
Note that in the situation above, $\Vec_L\cong\hbox{Rep}(L^*)$, where $L^*=\Hom(L,\Qlcl^*)$ is the Pontryagin dual. The automorphism group $\Gamma'$ of the fiber functor $\E\to \Vec$ is an extension of groups $0\to L^*\to \Gamma'\to \Gamma\to 0$. In particular, in case both $\Gamma$, $A$ are $p$-groups, it follows that $\Gamma'$ is one too.
\erk

\section{Appendix: The class $\mathfrak{C}_p^{\pm}$ of modular categories}\label{B}

In this appendix, we define and study some properties of the class $\mathfrak{C}_p^{\pm}$ of modular categories. We will follow \cite{DGNO2}, but we give an alternative equivalent definition. In the following definition, we use the notions of metric groups and Witt-equivalence that we discussed in Appendix \ref{mgassssbe} above. We also use the notions of central braided crossed categories and equivariantization, as well as that of the modular category $\M(A,\theta)$ associated to metric group $(A,\theta)$ defined in Appendix {\ref{mcassmg}} above.

\bdefn\label{defcppm}
\bit
\item[(i)] We say that a modular category $\C$ belongs to the class $\mathfrak{C}_p^{+}$ if $\C\cong \D^\Gamma$ as a modular category, where $\D$ is a (faithfully graded and positive spherical) central braided $\Gamma$-crossed category with identity component $\M(A,\theta)$ such that $\Gamma$ is a finite $p$-group and $(A,\theta)$ is a metric $p$-group Witt-equivalent to the trivial metric group.
\item[(ii)] We say that a modular category $\C$ belongs to the class $\mathfrak{C}_p^{-}$ if $\C\cong \D^\Gamma$ as a modular category, where $\D$ is a (faithfully graded and positive spherical) central braided $\Gamma$-crossed category with identity component $\M(A,\theta)$ such that $\Gamma$ is a finite $p$-group and $(A,\theta)$ is a metric $p$-group Witt-equivalent to the anisotropic metric group $(\F_{p^2},\zeta^{Nm(\cdot)})$.
\item[(iii)] The class $\mathfrak{C}_p^{\pm}$ is defined to be the union $\mathfrak{C}_p^{+}\cup \mathfrak{C}_p^{-}$.
\eit
\edefn

By Proposition \ref{core} we have the following:
\bprop\label{dgno2}
(i) A modular category $\C\in\mathfrak{C}_p^{+}$ if and only if it can be realized as the center  of a pointed fusion category (equipped with the positive spherical structure) whose group of isomorphism classes of simple objects is a finite $p$-group.\\
(ii) A modular category $\C\in \mathfrak{C}_p^{-}$ if and only if it can be realized as $\D^{\Gamma}$ where $\Gamma$ is a $p$-group and $\D$ is a faithfully graded central braided $\Gamma$-crossed category (equipped with positive spherical structure) with identity component $\M_p^{anis}=\M(\F_{p^2},\zeta^{Nm(\cdot)})$.
\eprop
\bpf
A pointed fusion category (with positive spherical structure) can be equivalently  thought of as a faithfully graded (positive spherical) central braided crossed  category with identity component $\Vec$ and that the equivariantization of such a category is equivalent to its Drinfeld center. (See \cite[\S 4.4.10]{DGNO}.) Both the statements now follow from Proposition \ref{core} and Remark \ref{pgroup}.
\epf

We now state another equivalent characterization of the class $\mathfrak{C}_p^{\pm}$ of modular categories proved in \cite{DGNO2}. We refer to \cite[\S2.4.1-3 and \S6]{DGNO} for the notions of Frobenius-Perron dimension, categorical dimension, Gauss sums and multiplicative central charge.

\bprop(\cite{DGNO2})
A modular category $\C$ belongs to $\mathfrak{C}_p^{\pm}$ if and only if the following three conditions are satisfied:
\bit
\item[(i)] The Frobenius-Perron dimension of $\C$ equals $p^{2k}$ for some $k\in \f{Z}^+$.
\item[(ii)] The categorical dimensions of all simple objects of $\C$ are positive integers.
\item[(iii)] The multiplicative central charge of $\C$ is $\pm 1$.
\eit
If $\C$ satisfies the conditions above then the multiplicative central charge of $\C$ is $1$ iff $\C\in \mathfrak{C}_p^{+}$ and the multiplicative central charge of $\C$ is $-1$ iff $\C\in \mathfrak{C}_p^{-}$.
\eprop

\bpf
Suppose that $\C\in\mathfrak{C}_p^{\pm}$. Hence by Proposition \ref{dgno2}, we have $\C\cong \D^\Gamma$ where $\D$ is some faithfully graded positive spherical central braided $\Gamma$-crossed category with trivial component  $\M(A,\theta)$ where the metric group $(A,\theta)$ is either $(0,1)$ (in case $\C\in\mathfrak{C}_p^{+}$) or $(\F_{p^2},\zeta^{Nm(\cdot)})$ (in case $\C\in\mathfrak{C}_p^{-}$) and $\Gamma$ is a finite $p$-group. The proof is very similar to the proofs of Propositions \ref{prop1} and \ref{cdpi}. We have that $\FPdim(\C)=\FPdim(\D^\Gamma)=|\Gamma|\cdot \FPdim(\D)=|\Gamma|^2\cdot\FPdim(\M(A,\theta))=|\Gamma|^2\cdot|A|$ which is an even power of $p$ since $|A|$ is either 1 or $p^2$ and $\Gamma$ is a finite $p$-group. The categorical dimensions of all simple objects of $\D$ are 1 since $\D$ is pointed and equipped with the positive spherical structure. Hence it follows that the categorical dimensions of all objects of $\D^\Gamma\cong \C$ must be positive integers. Finally, by \cite[Thm. 6.16]{DGNO} the Gauss sums of $\C$ are given by $\tau^{\pm}(\C)=\tau^{\pm}(\D^\Gamma)=|\Gamma|\cdot\tau^{\pm}(\M(A,\theta))$

\[
  =|\Gamma|\cdot \tau^{\pm}(A,\theta)=
  \begin{cases}
   |\Gamma| & \text{if } \C\in\mathfrak{C}_p^{+}\\
   -p\cdot|\Gamma|       & \text{if } \C\in\mathfrak{C}_p^{-}.
  \end{cases}
\]

Hence we see that the multiplicative central charge of $\C$ is 1 if $ \C\in\mathfrak{C}_p^{+}$ and -1 if $ \C\in\mathfrak{C}_p^{-}$.

\brk
The last statement implies that the classes $\mathfrak{C}_p^{+}$ and $\mathfrak{C}_p^{-}$ are disjoint.
\erk

We have completed the proof of the {\it ``only if''} part of the proposition. We refer to \cite{DGNO2} for the {\it ``if''} part.

\epf

The class $\mathfrak{C}_p^{-}$ may also be characterized as follows:
\bprop
(\cite{DGNO2})
A modular category $\C\in\mathfrak{C}_p^{-}$ if and only if the modular category $\C\boxtimes\M_p^{anis} \in \mathfrak{C}_p^{+}$.
\eprop

\bibliographystyle{ams-alpha}

\begin{thebibliography}{ABC}
\bibitem[B\'e]{Be} L. B\'egueri. {\it Dualit\'e sur un corps local \`a corps r\'esiduel alg\'ebriquement clos}, M\'em. Soc. Math. France (N.S.), 1980/81, no. 4.
\bibitem[Bo]{B} M. Boyarchenko. {\it Characters of Unipotent Groups over Finite Fields}, Selecta Mathematica, Vol. 16 (2010), No. 4, pp. 857--933, 	arXiv:0712.2614v4.
\bibitem[BD06]{BD06} M. Boyarchenko and V. Drinfeld. {\it A motivated introduction to character sheaves on unipotent groups in positive characteristic}, September 2006, 		arXiv:math/0609769v2.
\bibitem[BD08]{BD08} M. Boyarchenko and V. Drinfeld. {\it Character Sheaves on Unipotent Groups in Positive Characteristic: Foundations}, October 2008, 	arXiv:0810.0794v1.
\bibitem[Da] {Da} S. Datta. {\it Metric Groups attached to Skew-symmetric Biextensions}, Transformation Groups (2010), 	arXiv:0809.5082v2.
\bibitem[De] {De} T. Deshpande. {\it Heisenberg Idempotents on Unipotent Groups}, Math. Res. Lett. 17 (2010), no. 3, 415 - 434, arxiv:0907.3344.
\bibitem[DGNO]{DGNO} V. Drinfeld, S. Gelaki, D. Nikshych and V. Ostrik. {\it On Braided Fusion Categories I}, June 2009, 		arXiv:0906.0620v1. 
\bibitem[DGNO2]{DGNO2} V. Drinfeld, S. Gelaki, D. Nikshych and V. Ostrik. {\it On Braided Fusion Categories II}, In preparation.
\bibitem[EM]{EM} S. Eilenberg, S. MacLane. {\it Cohomology Theory in Abstract Groups II: Group Extensions with a non-Abelian Kernel}, Ann. of Math. Vol. 48, No.2 (1947), pp. 326-341.
\bibitem[ENO02]{ENO02} P. Etingof, D. Nikshych, V. Ostrik. {\it On Fusion Categories}, March 2002,  arXiv:math/0203060v10 [math.QA].
\bibitem[ENO09]{ENO} P. Etingof, D. Nikshych, V. Ostrik. {\it Fusion Categories and Homotopy Theory}, September 2009, arXiv:0909.3140v2.
\bibitem[H]{H} J. Humphreys. {\it Linear Algebraic Groups}, Springer-Verlag.
\end{thebibliography}

\end{document}